\documentclass[12pt]{amsart}
\usepackage{verbatim}
\usepackage{graphicx}

\title{Constant Mean Curvature Trinoids}

\author{N. Schmitt}
\address{Nicholas Schmitt, Center for Geometry, Analysis, Numerics and Graphics,
University of Massachusetts, Amherst, MA, USA}
\email{nick@gang.umass.edu}

\numberwithin{equation}{section}

\newcommand{\transpose}[1]{{{#1}^t}}
\newcommand{\interior}[1]{{#1^{\circ}}}

\newcommand{\abs}[1]{{\lvert#1\rvert}}
\newcommand{\norm}[1]{\lVert#1\rVert}

\newcommand{\half}{\tfrac{1}{2}}
\newcommand{\fourth}{\tfrac{1}{4}}

\newcommand{\ol}[1]{\overline{#1}}

\DeclareMathOperator{\Ad}{Ad}

\DeclareMathOperator{\Order}{O}

\DeclareMathOperator{\Span}{span}

\DeclareMathOperator{\ad}{ad}

\DeclareMathOperator{\cross}{\times}

\DeclareMathOperator{\del}{\partial\!}

\DeclareMathOperator{\id}{I}

\DeclareMathOperator*{\ord}{ord}

\DeclareMathOperator*{\res}{res}
\DeclareMathOperator{\re}{Re}

\DeclareMathOperator{\suchthat}{|}

\DeclareMathOperator{\tensor}{\otimes}
\DeclareMathOperator{\tr}{tr}
\DeclareMathOperator{\trace}{tr}
\DeclareMathOperator{\notequiv}{\equiv\!\!\!\!\!\slash\,}

\newcommand{\MatrixGroup}[1]{{\rm{#1}}}
\newcommand{\matGL}{\MatrixGroup{GL}}

\newcommand{\matSL}{\MatrixGroup{SL}}

\newcommand{\matSU}{\MatrixGroup{SU}}
\newcommand{\matU}{\MatrixGroup{U}}
\newcommand{\matgl}{\MatrixGroup{gl}}
\newcommand{\matsl}{\MatrixGroup{sl}}
\newcommand{\matsu}{\MatrixGroup{su}}
\newcommand{\mattwo}{\MatrixGroup{M}_{2\times 2}}

\theoremstyle{plain}
\newtheorem{theorem}{Theorem}[section]
\newtheorem{corollary}[theorem]{Corollary}
\newtheorem{lemma}[theorem]{Lemma}

\newtheorem*{corollary*}{Corollary}
\newtheorem*{lemma*}{Lemma}
\newtheorem*{proposition*}{Proposition}
\newtheorem*{theorem*}{Theorem}
\theoremstyle{definition}

\newtheorem{definition}[theorem]{Definition}

\newtheorem{notation}[theorem]{Notation}

\newtheorem{remark}{Remark}[section]
\newtheorem*{algorithm*}{Algorithm}
\newtheorem*{application*}{Application}
\newtheorem*{assertion*}{Assertion}
\newtheorem*{assumption*}{Assumption}
\newtheorem*{axiom*}{Axiom}
\newtheorem*{claim*}{Claim}
\newtheorem*{conjecture*}{Conjecture}
\newtheorem*{definition*}{Definition}
\newtheorem*{example*}{Example}
\newtheorem*{notation*}{Notation}
\newtheorem*{note*}{Note}
\newtheorem*{observation*}{Observation}
\newtheorem*{question*}{Question}
\newtheorem*{remark*}{Remark}
\theoremstyle{plain}

\newcommand{\bbC}{\mathbb{C}}

\newcommand{\bbP}{\mathbb{P}}

\newcommand{\bbR}{\mathbb{R}}
\newcommand{\bbS}{\mathbb{S}}

\newcommand{\bbZ}{\mathbb{Z}}

\newcommand{\calA}{\mathcal{A}}

\newcommand{\calC}{\mathcal{C}}
\newcommand{\calD}{\mathcal{D}}

\newcommand{\calM}{\mathcal{M}}

\newcommand{\calT}{\mathcal{T}}

\newcommand{\LoopG}[1]{\Lambda_{#1} G}
\newcommand{\LoopGL}[1]{\Lambda_{#1}\matGL_2(\bbC)}
\newcommand{\LoopSL}[1]{\Lambda_{#1}\matSL_2(\bbC)}

\newcommand{\LoopuG}[1]{\Lambda_{#1}^{*}G}
\newcommand{\LoopuGL}[1]{\Lambda_{#1}^{*}\matGL_2(\bbC)}
\newcommand{\LoopuSL}[1]{\Lambda_{#1}^{*}\matSL_2(\bbC)}

\newcommand{\LooppG}[1]{\Lambda_{#1}^{+}G}
\newcommand{\LooppGL}[1]{\Lambda_{#1}^{+}\matGL_2(\bbC)}
\newcommand{\LooppSL}[1]{\Lambda_{#1}^{+}\matSL_2(\bbC)}

\newcommand{\LoopprG}[1]{\Lambda_{#1}^{+,\bbR}G}
\newcommand{\LoopprGL}[1]{\Lambda_{#1}^{+,\bbR}\matGL_2(\bbC)}

\newcommand{\TriG}{\calT G}

\newcommand{\TriGL}{\calT GL_2(\bbC)}

\newcommand{\TrirG}{\calT^{\bbR}G}

\newcommand{\TrirGL}{\calT^{\bbR}GL_2(\bbC)}

\newcommand{\LoopprX}{\Lambda_{\uparrow 1}^{+,\bbR}\mattwo(\bbC)}

\newcommand{\Loopsl}[1]{\Lambda^{-1}_{#1}\matsl_2(\bbC)}
\newcommand{\Loopgl}[1]{\Lambda^{-1}_{#1}\matgl_2(\bbC)}
\newcommand{\Potentialsl}[2]{\Omega^1_{#1}(\Loopsl{#2})}
\newcommand{\Potentialgl}[2]{\Omega^1_{#1}(\Loopgl{#2})}

\newcommand{\gauge}[2]{{{#1}{.}{#2}}}

\DeclareMathOperator{\tracefree}{\mathsf{tracefree}}

\newcommand{\Sym}[2]{\mathsf{Sym}_{#1}[#2]}
\newcommand{\Uni}[2]{\mathsf{Uni}_{#1}[#2]}
\newcommand{\Pos}[2]{\mathsf{Pos}_{#1}[#2]}
\newcommand{\Note}[1]{}

\newcommand{\Wppp}{\mbox{\scriptsize$[+\!+\!+]$}}
\newcommand{\Wppm}{\mbox{\scriptsize$[+\!+\!-]$}}
\newcommand{\Wpmm}{\mbox{\scriptsize$[+\!-\!-]$}}
\newcommand{\Wmmm}{\mbox{\scriptsize$[-\!-\!-]$}}
\newcommand{\Del}{{\mbox{\rm \scriptsize 0}}}

\newcommand{\msmall}[1]{\text{\footnotesize $\displaystyle #1$}}

\begin{document}

\thanks{
{\indent This research was supported by National Science Foundation
grant DMS-00-76085.}
}

\begin{abstract}
This paper constructs a family of constant mean curvature immersions of the
thrice-punctured Riemann sphere into $\bbR^3$ with asymptotically
Delaunay ends via loop group methods.
\end{abstract}

\maketitle


\typeout{== intro.tex =============================================}

\section*{Introduction}\label{sec:intro}

A \emph{trinoid} is a conformal immersion of the
thrice-punctured Riemann sphere into $\bbR^3$ with constant mean
curvature (CMC) with asymptotically Delaunay ends.

CMC trinoids were first constructed in~\cite{kap1}.
The family of Alexandrov-embedded trinoids was classified
in~\cite{gks2}.

In this paper, a family of trinoids
is constructed via the Dorfmeister-Pedit-Wu (DPW)~\cite{dpw}
construction.
This family is three-dimensional,
parametrized by the asymptotic necksizes of the ends, and has four
connected components, according as the ends are embedded unduloids or
immersed nodoids.

Via the DPW construction, every CMC immersion can be obtained by first
solving a linear meromorphic ODE
\[
d\Phi_\lambda = \Phi_\lambda\xi_\lambda,\quad
\Phi_{\lambda}(z_0)=\Phi_\lambda^0.
\]
in 2-by-2 matrices which depend on loop parameter
$\lambda\in\bbS^1$.  A loop group factorization is then
applied to $\Phi=FB$ to produce the $\matSU_2$ frame $F$ for an
associate family of CMC immersions.

The first step in the DPW construction of CMC surfaces is to write down a
suitable family of potentials $\xi_\lambda$ for the ODE. To produce
asymptotically Delaunay ends, a natural choice is a potential which
is gauge equivalent to a linear superposition of three potentials for
Delaunay ends
\[
\xi_\lambda = \msmall{\sum} \frac{A_k(\lambda)}{z-z_k}dz.
\]
At each end, such a potential $\xi_\lambda$ has a simple pole whose residue
$A_k(\lambda)$
encodes the asymptotic necksize.
That the ends are asymptotic to half Delaunay surfaces follows from the
fact that at each end the potential is a perturbation of a potential
which produces a Delaunay surface.

Since the solution $\Phi_\lambda$ to the ODE has
monodromy around the ends, it is necessary to simultaneously close the
ends by choosing a suitable initial condition $\Phi_\lambda^0$. The
ends are closed when the monodromy representation of $\Phi_\lambda$ is
unitary. In the case of three ends, the necessary and sufficient
condition for unitarizing the monodromy representation is the triangle
inequalities on the 2-sphere
\begin{align*}
  & \nu_1+\nu_2+\nu_3\le 1\\
  & \nu_i\le \nu_j+\nu_k
\end{align*}
where $\nu_1,\,\nu_2,\,\nu_3$ are the necksizes of the three ends,
which can be read off from the resides of $\xi_\lambda$.  (In the case
of $n>3$ ends, the spherical $n$-gon inequalities on the necksizes
are a necessary condition.)  The
$\Phi_\lambda^0$ which unitarizes the monodromy pointwise in $\lambda$
is then ``glued'' holomorphically in $\lambda$
(theorem~\ref{thm:glue}) to produce an initial condition for which the
ends of the CMC immersion are closed.

The tools developed here
will be useful for constructing further examples of CMC surfaces,
including $n$-noids in $\bbR^3$,
$n$-noids in $H^3$ and $\bbS^3$, and $n$-noids dressed by
B\"ackland transformations.

The images in this paper were generated with \texttt{CMCLab}, a
numerical implementation of the DPW algorithm developed by the author
using algorithms for loop group factorizations explicated by
I.~McIntosh.  The software, documentation, and a gallery of CMC
surface images produced by the software is available at the Center for
Geometry, Analysis, Numerics and Graphics (GANG) website,
\texttt{www.gang.umass.edu}.

\textbf{Acknowledgments.}
I am grateful to F.~Pedit, I.~McIntosh, R.~Kusner, W.~Rossman, and
M.~Kilian for helpful discussions, and to J.~Dorfmeister and
H.~Wu for making available their work on trinoids in preprint
form~\cite{dw-trinoid}.

\subsection*{Outline of the paper}
\mbox{}

Section~\ref{sec:prelim} (Preliminary) explicates the DPW construction
and provides background theory relating to monodromy, closing conditions
and gauge equivalence.

Section~\ref{sec:perturb} (Delaunay immersions) 
gives the basic background concerning Delaunay surfaces,
their dressings and their asymptotic growth rates.

Section~\ref{sec:asymptotics} (Perturbations of Delaunay immersions)
proves the asymptotics theorem~\ref{thm:asymptotic}, showing
that the CMC immersion arising from a perturbation of a Delaunay
DPW potential has an asymptotically Delaunay end.

Section~\ref{sec:glue} (Unitarization of loop group
monodromy representations) proves the gluing theorem~\ref{thm:glue},
that under the assumption of the pointwise unitarizability of a monodromy
representation on $\bbS^1_\lambda$, there exist an $r$-dressing which
conjugates the monodromy representation to an $r$-unitary representation.

In section~\ref{sec:trinoid} (Constructing trinoids)
a family of trinoid potentials $\calT$ is given
(definition~\ref{def:trinoid}).
Theorem ~\ref{thm:unitary3} gives a necessary and sufficient
condition for simultaneous unitarizability of a monodromy representation
on the thrice-punctured sphere in terms of the eigenvalues.
Theorem~\ref{thm:triunitary}
shows that the monodromy representation for $\xi\in\calT$
is unitarizable pointwise on $\bbS^1_\lambda$.
Theorems~\ref{thm:main} and~\ref{thm:symmetry} draw on these results to
construct a family of trinoids parametrized by the three end weights.



\begin{figure}[ht]
  \centering
  \includegraphics[width=110pt]{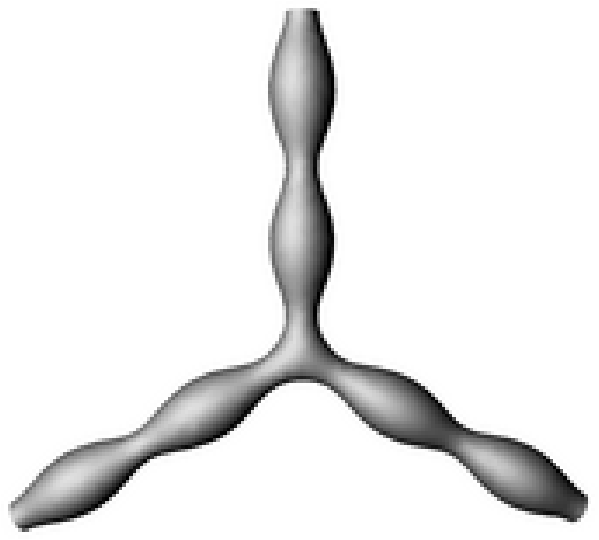}
  \includegraphics[width=110pt]{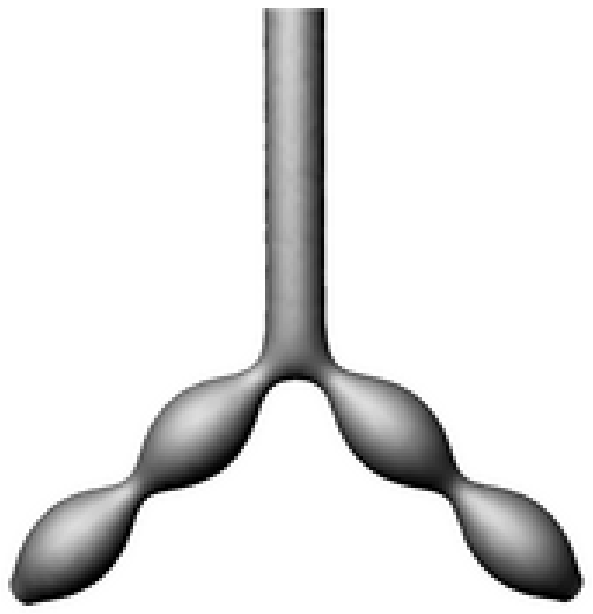}
  \includegraphics[width=110pt]{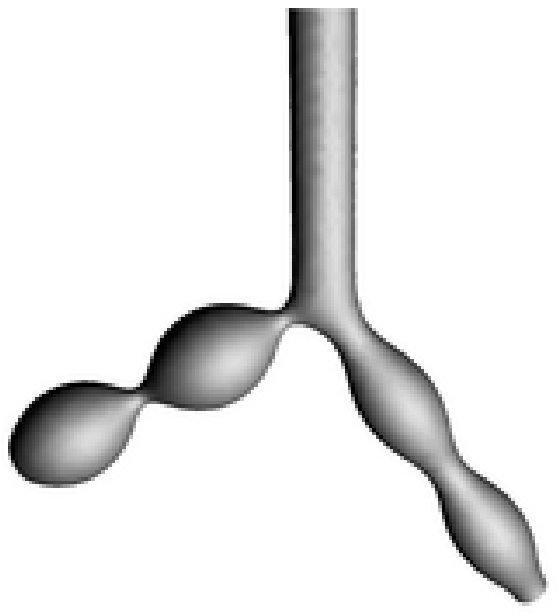}
\caption{Equilateral, isosceles and scalene CMC trinoids with unduloid ends.
Their respective necksizes are
($\tfrac{1}{3},\,\tfrac{1}{3},\,\tfrac{1}{3}\bigr)$,
($\tfrac{1}{2},\,\tfrac{1}{4},\,\tfrac{1}{4}\bigr)$
and
($\tfrac{1}{2},\,\tfrac{1}{3},\,\tfrac{1}{6}\bigr)$
These examples are static in the sense that the sum of their necksizes is
the maximum of 1.}
\end{figure}

\typeout{== prelim.tex =============================================}
\section{Preliminary}\label{sec:prelim}

\subsection{The DPW construction}
\label{sec:dpw}

\begin{notation}
  The following notation is used for circles, disks and annuli
  in the domain $\bbP^1_\lambda$ of the loop parameter $\lambda$.
  Let $r\in (0,1]$.
  \begin{align*}
    \calC_r &= \{\lambda\in\bbC\suchthat\abs{\lambda}=r\}\\
    \calD_r &= \{\lambda\in\bbC\suchthat\abs{\lambda}<r\}\\
    \calD_r' &= \{\lambda\in\bbC\suchthat\abs{\lambda}>r\}\cup\{\infty\}\\
    \calA_r &= \{\lambda\in\bbC\suchthat r<\abs{\lambda}<1/r\}.
  \end{align*}
\end{notation}

\begin{notation}
  For a map $X:A_r\to M_{k\times k}(\bbC)$, the star operator
  is defined as
  \begin{equation}\label{eq:star}
    X^{*}(\lambda) = \transpose{\overline{X(\overline{\lambda}^{-1})}}.
  \end{equation}
\end{notation}

\begin{notation}
  The following groups are defined for $G$ either $\matGL_n(\bbC)$ or
  $\matSL_n(\bbC)$.

  Let $\TriG$ denote the group of upper triangular elements of $G$
  and $\TrirG$ the subgroup of $\TriG$ whose diagonal elements
  are in $\bbR^{+}$.
  For $r\in(0,\,1]$,

  $\bullet$
  $\LoopG{r}$ (``loops'') is the group of analytic maps $C_r\to G$.

  $\bullet$
  $\LoopuG{1}$ (``unitary loops'') is the subgroup of $\LoopG{1}$ of loops
  $X$ which satisfy the reality condition
  \begin{equation}\label{eq:reality}
    \widehat{X}={{\left.\widehat{X}\right.}^\ast}^{-1}.
  \end{equation}
  For $r\in(0,\,1)$,
  $\LoopuG{r}$ (``$r$-unitary loops'')
  is the subgroup of $\LoopG{r}$ of loops $X$
  such that $X$ is the boundary of an analytic map
  $\widehat{X}:\calA_r\to G$ satisfying
  the reality condition.

  $\bullet$
  $\LooppG{r}$ is the subgroup of $\LoopG{r}$ of loops $X$ such that
  $X$ is the boundary of an analytic map $\widehat{X}:\calD_r\to G$
  satisfying $\widehat{X}(0)\in \TriG$.

  $\bullet$
  $\LoopprG{r}$ (``positive loops'')
  is the subgroup of $\LooppG{r}$ of loops $X$ such that $X$ is the
  of boundary of an analytic map $\widehat{X}:\calD_r\to G$ satisfying
  $\widehat{X}(0)\in T^\bbR G$.

  $\bullet$
  $\LoopprX$ is the set of analytic maps
  $X:\bbS^1\to\mattwo(\bbC)$ such that $X$ is the boundary of an
  analytic map $\widetilde{X}:\calD_1\to\matGL_2(\bbC)$ satisfying
  $\widehat{X}(0)\in \TrirGL$.
\end{notation}

The $r$-Iwasawa factorization~\cite{PressleySegal, McI2} is as follows.

\begin{theorem*}[Iwasawa factorization theorem]
  Let $r\in(0,\,1]$. Take $G$ to be either $\matGL_2(\bbC)$ or
  $\matSL_2(\bbC)$.
  Then any $X\in\LoopG{r}$ can be factored uniquely
  \[
  X = X_u X_{+}
  \]
  with $X_u\in\LoopuG{r}$ and $X_{+}\in\LoopprG{r}$.
  The induced map
  \[\LoopG{r}\to\LoopuG{r}\cross\LoopprG{r}\]
  is an analytic diffeomorphism
  \Note{with respect to what structure}.
\end{theorem*}

The projections of the $r$-Iwasawa factorization of $X$ to
the first and second factors are respectively denoted
by $\Uni{r}{X}$ and $\Pos{r}{X}$.
For loops $F\in\LoopG{r}$ and $C\in\LoopG{r}$, the $r$-dressing action
of $C$ on $F$ is $\Uni{r}{CF}$.

\begin{notation}
  $\Loopgl{r}$ and $\Loopsl{r}$ are respectively the sets of holomorphic
  $\matgl_2(\bbC)$- and $\matsl_2(\bbC)$-valued
  functions on $\calD_r^\ast$ which extend meromorphically to
  $\lambda=0$ and whose expansion in $\lambda$ at $\lambda=0$
  is of the form
  \[
  \begin{pmatrix}0 & \alpha\\0 & 0\end{pmatrix}\lambda^{-1}+
  \Order(\lambda^0).
  \]

  For a Riemann surface $\Sigma$ and complex vector space $V$,
  $\Omega^{1}_\Sigma(V)$ denotes the holomorphic $V$-valued
  1-forms on $\Sigma$.
\end{notation}

For $X\in\LoopG{r}$, the notation $X'$ means differentiation with respect
to $\theta$, where $\lambda=e^{i\theta}$.
We have ${(X')}^\ast={(X^\ast)}'$.

For $X\in\matgl_2(\bbC)$, $\tracefree(X)=X-\half(\trace{X})\id$.

The DPW construction~\cite{dpw} is as follows.

\begin{theorem*}[DPW]
  Let $\Sigma$ be a Riemann surface and $\widetilde{\Sigma}$ its universal
  cover. Let $r\in(0,\,1]$.
  Let $\xi\in\Potentialgl{\Sigma}{r}$.
  Let $z_0\in \widetilde{\Sigma}$ and let $\Phi_0\in\LoopGL{r}$.
  Let $\Phi:\widetilde{\Sigma}\to\LoopGL{r}$ be the solution
  to the initial value problem
  \begin{equation}\label{eq:ode}
    d\Phi = \Phi \xi;\quad \Phi(z_0) = \Phi_0.
  \end{equation}
  This initial value problem is denoted by the triple $(\xi,\,z_0,\,\Phi_0)$.

  Let
  \begin{equation*}
    \Phi = F B
  \end{equation*}
  be the $r$-Iwasawa factorization of $\Phi$.
  Then $F$ extends to a map $F:\widetilde{\Sigma}\to\LoopuGL{r}$
  and $\left.F\right|_{\bbS^1}$ takes values in $\matU_2$.
  $F$ is called the \emph{extended frame}.

  Let $\Sym{\lambda}{\,\cdot\,}$ be defined on maps
  $F:\widetilde{\Sigma}\to\LoopuGL{1}$ by
  \begin{equation}\label{eq:sym}
    \Sym{\lambda}{F} = -2H^{-1}\tracefree\bigl(F'F^{-1}\bigr).
  \end{equation}
  For each $\lambda\in\bbS^1$, the map $\Sym{\lambda}{F}$ is
  a conformal constant mean curvature immersion
  $\Sigma\to\matsu_2\equiv\bbR^3$
  with mean curvature $H$.
  The family $\Sym{\lambda}{F}$ over $\lambda\in\bbS^1$
  is an \emph{associate family} of CMC immersions.
\end{theorem*}

\begin{remark}
  \label{rem:dpw}
The Hopf differential of $f_\lambda$ is
$-2H^{-1}\alpha\beta\lambda^{-1}$ and
its metric is $4H^{-2}R^2\alpha\tensor\ol{\alpha}$,
where $R=B_{11}/B_{22}$ and $(B_{ij})= B|_{\lambda=0}$.
\end{remark}

\subsection{Monodromy and closing conditions}
\label{sec:close}

\begin{lemma}
  \label{thm:rindependent}
  Let $0<r_1<r_2\le1$, and suppose that $\Phi_1\in\LoopGL{r_1}$
  $\Phi_2\in\LoopGL{r_2}$ are the boundary of an analytic map
  $\Phi:\{r_1<\abs{\lambda}<r_2\}$.
  Let $\Phi_j=F_jB_j$ be the $r_j$-Iwasawa factorizations
  of $\Phi_j$, $j=1,\,2$.
  Let $F_j$ be the extension of $F_j$ to $\calA_{r_j}$
  and $B_j$ be the extension of $B_j$ to $\calD_{r_j}$.
  Then
  $F_2$ extends analytically to $\calA_{r_1}$ and is equal to
  $F_1$ there, and
  $B_1$ extends analytically to $\calD_{r_2}$ and is equal to
  $B_2$ there.
\end{lemma}

\begin{proof}
  Since $\Phi$ and $B_2$ are analytic on $\{r_1<\abs{\lambda}<r_2\}$,
  $F_2=\Phi_2 B_2^{-1}$ extends analytically to $\calA_{r_1}$.
  Since $\Phi$ and $F_1$ are analytic on $\{r_1<\abs{\lambda}<r_2\}$,
  $B_1=F_1^{-1}\Phi_1$ extends analytically to $\calD_{r_2}$.
  Hence $\Phi=F_1B_1=F_2B_2$ is an $r$-Iwasawa factorization for any
  $r\in[r_1,\,r_2]$, so by the uniqueness of $r$-Iwasawa factorization,
  $F_1=F_2$ and $B_1=B_2$.
\end{proof}

\begin{notation}\label{not:monodromy}
  Let $\Sigma$ be a Riemann surface, $\widetilde{\Sigma}\to\Sigma$
  its universal cover, and  $\Gamma$ the group of deck transformations
  for this cover.
  Let $r\in(0,\,1]$, let $\xi\in\Potentialgl{\Sigma}{r}$ and
  let $\Phi:\widetilde{\Sigma}\to\LoopGL{r}$
  be a solution to the ODE $d\Phi=\Phi\xi$.
  The \emph{monodromy representation} of $\Phi$ is the map
  $M_\Phi:\Gamma\to\LoopGL{r}$ defined by
  $M_\Phi(\tau)=(\tau^\ast\Phi)\Phi^{-1}$.

  In the case  $\Sigma=\Sigma_0\setminus\{p_1,\dots,p_n\}$ for a
  closed Riemann surface $\Sigma_0$,
  we define the ``monodromy of $\Phi$ at $p_k$'' as $M_\Phi(\tau)$,
  where $\tau\in\Gamma$ is defined as follows:
  let $U$ be an annular neighborhood of $p_k$,
  $\gamma:[0,\,1]\to U$ a closed curve with winding number $1$ around $p_k$,
  and $\tau\in\Gamma$ the deck transformation satisfying
  $\tau(\gamma(0))=\gamma(1)$.
\end{notation}

\begin{lemma}\label{thm:close}
  Let $\Sigma$, $\widetilde{\Sigma}$, $\Gamma$,
  $r$, $\xi$, $\Phi$ and $M$ be as
  in notation~\ref{not:monodromy}, and suppose that
  \begin{equation}
    \label{eq:close1}
    M\in\LoopuGL{r}.
  \end{equation}
  Let $\lambda_0\in\bbS^1$,
  and let $f_{\lambda}=\Sym{\lambda}{\Uni{r}{\Phi}}$.
  Let $\tau\in\Gamma$.
  Then $(\tau^\ast F)F^{-1}$ is $z$-independent,
  $(\tau^{*}B) B^{-1}=\id$, and the following are equivalent:
  \begin{equation}
    \label{eq:close2}
    \begin{split}
      &M(\tau,\,\lambda_0) \text{ is a multiple of $\id$}\\
      &\tracefree(M'(\tau,\,\lambda_0)) = 0
    \end{split}
  \end{equation}
  and
  \begin{equation}
    \label{eq:fmono}
    \tau^{*}f_{\lambda_0}=f_{\lambda_0},
  \end{equation}
  where $f=\Sym{}{F}$.
  \Note{Need some non-degeneracy on $f$.}
  In the case $M\in\LoopuSL{r}$,
  conditions~\eqref{eq:close2} are equivalent to
  \begin{equation}
    \label{eq:close2a}
    M(\tau,\,\lambda_0) = \pm\id,\quad M'(\tau,\,\lambda_0) = 0.
  \end{equation}
\end{lemma}

\begin{proof}
  Let $\Phi=FB$ be the $r$-Iwasawa factorization of $\Phi$. Then
  \begin{equation}
    \label{eq:cl}
    (\tau^{*}F^{-1}) (M(\tau)) F = (\tau^{*}B) B^{-1}.
  \end{equation}
  holds on $\calC_r$.
  Since $M\in\LoopuGL{r}$, the left hand side of
  equation~\eqref{eq:cl}
  is in $\LoopuGL{r}$ while the
  right hand side is in $\LoopprGL{r}$. The uniqueness of the $r$-Iwasawa
  factorization implies that each side of the equation is $\id$, so
  $M(\tau)= (\tau^\ast F)F^{-1}$
  and $(\tau^{*}B) B^{-1}=\id$ on $\calC_r$.
  Under the assumption~\eqref{eq:close1},
  $M(\tau)$ is in $\LoopuGL{r}$, so it
  extends analytically to $\calA_r$.

  We have
  \begin{equation}
    \label{eq:tauf}
    \tau^{*}f = M(\tau)fM(\tau)^{-1} -
    2H^{-1}\tracefree(M(\tau)'M(\tau)^{-1}).
  \end{equation}
  Assuming equation~\eqref{eq:close2},
  $M(\tau,\,\lambda_0)=\pm\id$ and $M'(\tau,\,\lambda_0)=0$,
  so the formula~\eqref{eq:tauf} evaluated at $\lambda_0$ yields
  $\tau^{*}f_{\lambda_0}=f_{\lambda_0}$.

  Conversely, note that for fixed $\lambda\in\bbS^1$, 
  the action on $f_\lambda$ defined by the right hand side of
  equation~\eqref{eq:tauf} is an isometry of $\matsu_2$.
  If equation~\eqref{eq:fmono} holds,
  then this isometry fixes $f_{\lambda_0}$ pointwise, so either
  $f_{\lambda_0}$ lies in two-dimensional subspace of $\matsu_2$
  or the isometry is the identity.
  Equations~\eqref{eq:close2} follow.
\end{proof}

The following lemma shows that condition~\ref{eq:close2} can be replaced
by an analogous condition on the eigenvalues of $M_\Phi$.

\begin{lemma}\label{thm:monoeigen}
  Let $\gamma$ be an open segment of $\bbS^1$,
  $M:\gamma\to\matU_2$
  an analytic map,
  $\rho_1,\,\rho_2$ the eigenvalues of $M$,
  and $\lambda_0\in\gamma$. Then
  the conditions~\eqref{eq:close2} are equivalent to
  \begin{equation}
    \label{eq:eigenclose}
    \rho_1(\lambda_0)=\rho_2(\lambda_0),\quad
    \rho_1'(\lambda_0)=\rho_2'(\lambda_0).
  \end{equation}
  In the case $M:\gamma\to\matSU_2$, these are equivalent to
  \begin{equation}
    \label{eq:eigenclosea}
    \rho_1(\lambda_0)=\pm 1,\quad \rho_1'(\lambda_0)=0.
  \end{equation}
\end{lemma}

\begin{proof}
  Since $M(\lambda_0)\in\matU_2$,
  $M(\lambda_0)$ is a multiple of $\id$ iff
  $\rho_1(\lambda_0)=\rho_2(\lambda_0)$.
  Assuming this, differentiating the characteristic equation
  $\rho^2-(\tr M)\rho+\det M=0$ twice
  and evaluating at $\lambda_0$ yields
  \begin{equation}
    \label{eq:Mrho}
    \rho'^2(\lambda_0)-(\tr M'(\lambda_0))\rho'(\lambda_0)+
    \det M'(\lambda_0)=0.
  \end{equation}
  Hence $\rho_1'(\lambda_0),\,\rho_2'(\lambda_0)$ are the eigenvalues
  of $M'(\lambda_0)$.
  But if $\tracefree(M'(\lambda_0))=0$, then the eigenvalues of
  $M'(\lambda_0)$ are equal.

  Conversely,
  since the eigenvalues of $M'(\lambda_0)$ are
  $\rho_j(\lambda_0)$, we have by equation~\eqref{eq:Mrho}
  \begin{equation}
    \label{eq:Mrho2}
    M'(\lambda_0)^2 - (\tr M'(\lambda_0))M'(\lambda_0)+
    (\det M'(\lambda_0))M'(\lambda_0)=0.
  \end{equation}
  If $\rho_1'(\lambda_0)=\rho_2'(\lambda_0)$, then
  $(\tr M'(\lambda_0))^2 = 4\det M'(\lambda_0)$ and
  equation~\eqref{eq:Mrho2} becomes
  $(\tracefree(M'(\lambda_0)))^2=0$.
  Differentiating $MM^\ast=\id$ shows that
  $\rho_1(\lambda_0)^{-1}M'(\lambda_0)$ is skew-hermitian.
  It follows that $M'(\lambda_0)$ a multiple of $\id$.
\end{proof}

\subsection{Gauge equivalence}
\label{sec:gauge}

\begin{notation}
  Let $\Sigma$ be a Riemann surface,
  $\xi\in\Potentialgl{\Sigma}{r}$ and
  $g:\widetilde{\Sigma}\to\LoopGL{r}$
  and suppose that the monodromy group of $g$ is
  a subgroup of $\bbC^\ast\id$.
  The gauged potential $\gauge{\xi}{g}$ is
  \begin{equation*}
    \gauge{\xi}{g} = g^{-1}\xi g + g^{-1}dg.
  \end{equation*}
\end{notation}
\noindent If $\Phi$ is a solution to the ODE $d\Phi=\Phi\xi$, then
$\Psi=\Phi g$ is a solution to the gauged ODE $d\Psi=\Psi(\gauge{\xi}{g})$.

The following lemma provides the basic facts relating to gauge equivalence.

\begin{lemma}\label{thm:gauge}
  Let $\Sigma$, $\widetilde{\Sigma}$, $\Gamma$,
  $r$, $\xi$, $\Phi$ and $M_\Phi$ be as
  in notation~\ref{not:monodromy}.
  Let $g:\widetilde{\Sigma}\to\LooppGL{r}$ (resp. $\LooppGL{r}$)
  and suppose that the monodromy of $g$ takes values in $\bbC^\ast\id$.
  Let $M_{\Phi g}$ be the monodromy of $\Phi g$. Then

  (i) $\gauge{\xi}{g}\in\Potentialsl{\Sigma}{r}$
  (resp $\Potentialsl{\Sigma}{r}$).

  (ii) $M_{\Phi g}=c M_{\Phi}$, where
  $c$ is the monodromy of $g$.

  (iii)
  $
  \Sym{\lambda}{\Uni{r}{\Phi}}=\Sym{\lambda}{\Uni{r}{\Phi g}}.
  $
\end{lemma}

\begin{proof}
  To show (i), an examination of the series for $g$ and $\xi$ in $\lambda$
  at $\lambda=0$ show that $\gauge{\xi}{g}$ is holomorphic at $\lambda=0$,
  hence $\gauge{\xi}{g}$ is holomorphic in $\calD_r$.
  In the case $\Potentialsl{\Sigma}{r}$, note that if $\det g=\id$ and $\xi$
  is tracefree, then $\gauge{\xi}{g}$ is tracefree.

  Proof of (ii):
  \begin{equation*}
    \begin{split}
      M_{\Phi g}(\tau) &= (\tau^\ast(\Phi g)) (\Phi g)^{-1} =
      (\tau^\ast\Phi) ((\tau^\ast g) g^{-1}) \Phi^{-1}\\
      &= c(\tau)(\tau^\ast\Phi)\Phi^{-1} =
      c(\tau)M_\Phi(\tau).
    \end{split}
  \end{equation*}

  \emph{Proof of (iii).}
  Let $\Phi=FB$ be the $r$-Iwasawa factorization of $\Phi$.
  Let $(BG)(0)=UT$ be the pointwise Iwasawa factorization
  of $(BG)(0)$ (so $U\in\matU_2$ and $T\in\TrirGL$).
  Then the $r$-Iwasawa factorization of $\Phi g$ is
  \begin{equation*}
    \Phi g =(FU)(U^{-1}Bg),
  \end{equation*}
  and
  \begin{equation*}
    \begin{split}
      \Sym{\lambda}{\Uni{r}{\Phi g}} &= \Sym{\lambda}{FU}
      = -2H^{-1}(FU)'(FU)^{-1}\\
      &= -2H^{-1}F'F^{-1}
      = \Sym{\lambda}{F} = \Sym{\lambda}{\Uni{r}{\Phi}}.
    \end{split}
  \end{equation*}
\end{proof}

\begin{figure}[ht]
  \centering
  \includegraphics[width=175pt]{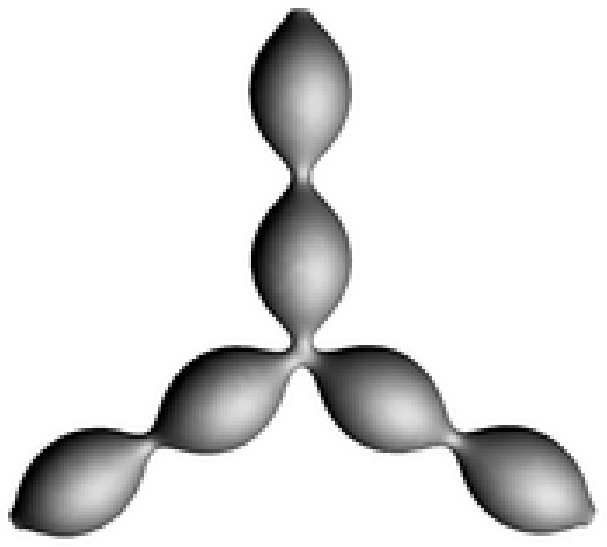}
  \includegraphics[width=175pt]{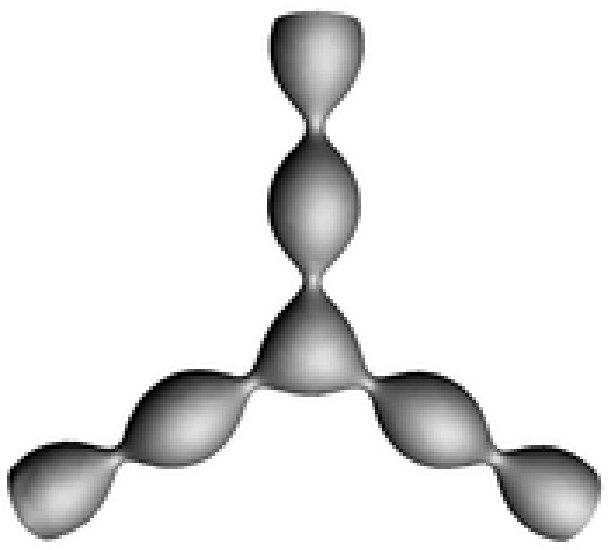}
\caption{A pair of CMC trinoids with unduloid ends (necksizes
$\bigl(\tfrac{1}{6},\,\tfrac{1}{6},\,\tfrac{1}{6}\bigr)$).
They have respectively a central neck and central bulge, exhibiting
a phase shift.}
\end{figure}

\typeout{== delaunay.tex =============================================}
\section{Delaunay immersions}
\label{sec:perturb}

\subsection{Delaunay surfaces via DPW}

CMC surfaces whose ends are asymptotic to Delaunay surfaces can be
constructed as local perturbations of a base Delaunay
surface. Hence Delaunay surfaces are first discussed.

The only CMC surfaces of revolution are the round cylinder,
the Delaunay unduloids (embedded Delaunay surfaces),
the sphere, and the Delaunay nodoids (immersed non-embedded Delaunay surfaces).

\begin{definition}
  Let $f$ be an conformal CMC immersion of constant mean curvature $H$
  which is a surface of revolution (a sphere, cylinder or Delaunay surface).
  The \emph{necksize} $n$ of $f$ is the minimum radius of the foliating
  circles, taken to be negative in the case of nodoids.
  The \emph{weight} of $f$ is $w = 4n(H^{-1}-n)$.
\end{definition}

In the case $H=1$,
the round cylinder has weight $1$ and necksize $\half$,
the unduloids have weight
in $(0,\,1)$ and necksize in $(0,\,\half)$,
the round sphere has weight and necksize $0$, and
the nodoids have weight and necksize in $(-\infty,\,0)$.

The DPW construction of Delaunay surfaces~\cite{kms} are as
follows.

\begin{lemma}\label{thm:del}
  Let $A\in\Loopsl{1}$ satisfy $A^\ast=A$, so that
  \begin{equation}
    \label{eq:del}
    A =
    \begin{pmatrix}
      c & a\lambda^{-1}+\ol{b}\\
      b+\ol{a}\lambda & -c
    \end{pmatrix}\quad
    a,\,b\in\bbC^\ast,\ c\in\bbR,
  \end{equation}
  Let $\Phi:\bbC\times\calA_0\to\matSL_2(\bbC)$
  be defined by $\Phi=\exp(\zeta A)$. Then

  (i)
  $f_\lambda=\Sym{\lambda}{\Uni{r}{\Phi}}$ is independent of the
  choice of $r\in(0,\,1]$.
  For each $\lambda\in\bbS^1$, $f_\lambda$ has screw symmetry.

  (ii)
  Let $\mu(\lambda)$ be an eigenvalue of $A$.
  If
  \begin{equation}
    \label{eq:del-close}
    \mu(1)=\pm\half,\quad
    \mu'(1)=0,
  \end{equation}
  then $f_{1}$ satisfies $f_{1}(\zeta+2\pi i)=f_{1}(\zeta)$
  and is a once-wrapped conformal immersion of a Delaunay surface with weight
  $16abH^{-2}$.
  The eigenvalues of the monodromy $M=\exp(2\pi i A)$ of $\Phi$
  are $\exp(\pm2\pi i \mu_w)$, where
  $w=16 ab\in(-\infty,\,1]\setminus\{0\}$ and
  \begin{equation}
    \label{eq:mu}
    \mu_w = \frac{1}{2}\sqrt{1+\frac{w(\lambda-1)^2}{4\lambda}}.
  \end{equation}

  (iv)
  If $A_1,\,A_2$ are of the form~\eqref{eq:del} with
  $\det A_1=\det A_2$, and
  $f_1=\Sym{}{\Uni{r}{\exp(z A_1}}$ and
  $f_2=\Sym{}{\Uni{r}{\Phi(z A_2}}$,
  then there exists an isometry $T:\matsu_2\to\matsu_2$
  and a coordinate change $z=\tilde{z}+c$, $c\in\bbR$
  such that $f_2(z)=T(f_1(\tilde{z}))$.
\end{lemma}

\begin{proof}
  $f$ is independent of the choice of $r$ by lemma~\ref{thm:rindependent}.
  Let $\theta\in\bbR$, $u=\exp(i\theta)$ and
  $U=\exp(i\theta A)$. Then $u\in\matSU_1$, $U\in\LoopuSL{r}$.
  Let $F=\Uni{r}{\Phi}$. Then
  \begin{equation*}
    \begin{split} 
      \Phi(uz) &= U\Phi(z)\\
      F(uz) &= UF(z)\\
      f(uz) &= \Ad_U f(z)-\tfrac{2}{H}U'U^{-1}.
    \end{split}
  \end{equation*}
  A calculation shows that there exists
  $S,\,T:\bbS^1\to\matsu_2$ such that
  \begin{equation*}
    f_1(uz) +S = \Ad_{U(1)}( f_1(z)+S ) + T.
  \end{equation*}
  This implies that $f$ has screw symmetry, and is hence an associate family
  of Delaunay immersions.
  The monodromy of this solution, $M_\Phi = \exp(2\pi i A)$,
  satisfies the closing condition~\eqref{eq:close1}.
  Under the hypotheses on the eigenvalues, $M_\Phi$ satisfies
  the closing condition~\eqref{eq:close2} at $\lambda_0=1$,
  so $f_1$ is monodromy-free along a loop around $z=0$.
  A calculation shows that the weight of $f_1$ is $16ab$.
  The proof of (iv) is omitted.
\end{proof}

\subsection{Dressed Delaunay immersions}
\label{sec:dress-del}

\begin{lemma}\label{thm:extend-unitary}
  Let $\gamma$ be an open segment of $\bbS^1$, $\lambda_0\in\gamma$,
  $\gamma^\ast=\gamma\setminus\{\lambda_0\}$
  and $M:\gamma^\ast\to\matU_2$ (resp. $\matSU_2$) a real analytic map which
  extends meromorphically to a neighborhood of $\gamma$.
  Then $M$ extends to a real analytic map $M:\gamma\to\matU_2$
  (resp. $\matSU_2$).
\end{lemma}

\begin{proof}
  Since $M$ takes values in $\matU_2$ on $\gamma^\ast$,
  its entries are bounded in absolute value by $1$ there.
  Since a meromorphic function
  at a pole is unbounded along every curve into the pole, the entries
  of $M$ cannot have poles at $\lambda_0$. Hence $M$ extends real
  analytically to $\lambda_0$.

  Since $MM^\ast=\id$ on $\gamma^\ast$,
  then $MM^\ast=\id$ on $\gamma$ by the continuity of $MM^\ast$.
  If $\det M(\lambda_0)=1$ on $\gamma^\ast$, then
  then $\det M=\id$ on $\gamma$ by the continuity of $\det M$.
\end{proof}

Lemma~\ref{thm:unitcomm} provides a ``unitary-commutator'' factorization
theorem, used in lemma~\ref{thm:dressdel}.

\begin{notation}
  \label{not:adjoint}
  For $X=\begin{pmatrix}x_{11} & x_{12} \\ x_{21} & x_{22}
  \end{pmatrix}\in\mattwo(\bbC)$,
  define $\widehat{X}=\begin{pmatrix}x_{22} &- x_{12} \\ -x_{21} & x_{11}
  \end{pmatrix}$, so $X\widehat{X}=(\det{X})\id$.
\end{notation}

\begin{lemma}\label{thm:unitcomm}
  Let $M\in\LoopuSL{1}$.
  Let $C:\bbS^1\to\mattwo(\bbC)$ be a real analytic map
  with $\det C\notequiv 0$
  such that the extension of $CMC^{-1}$ across $\{\det C=0\}$ is
  in $\LoopuSL{1}$ (see lemma~\ref{thm:extend-unitary}).
  Then there exists
  $U\in\LoopuSL{1}$ and $R:\bbS^1\to\mattwo(\bbC)$
  such that $C=UR$ and $[R,\,M]=0$.
\end{lemma}

\begin{proof}
  Since $C\notequiv 0$, there exists $c\in\bbC^\ast$ such that
  $V = c C+\widehat{c C}^\ast\notequiv 0$.
  Then $V=\widehat{V}^\ast$. It follows that
  $\det V$ takes values in $\bbR^{\ge 0}$ and
  $\det V\notequiv 0$.
  Hence there exists a well-defined non-negative square root $\sqrt{\det V}$
  on $\bbS^1$ which is not identically $0$.

  Define $U=(\det V)^{-1/2}V$ away from the zeros of $\det V$.
  By lemma~\ref{thm:extend-unitary}(i), $U$ extends analytically
  to $\bbS^1$ and $U\in\LoopuSL{1}$.

  Define $R=U^{-1}C$. Then
  $CMC^{-1}C-CM=0$ and, using the fact that $M$ and $CMC^{-1}$ satisfy
  the reality condition, $CMC^{-1}\widehat{C}^\ast-\widehat{C}^\ast M=0$.
  Hence $CMC^{-1}U-UM=0$, and so
  \[
  [R,\,M]=[U^{-1}C,\,M] = U^{-1}(CMC^{-1}U-UM)U^{-1}C = 0.
  \]
  Hence $U$, $R$ satisfy the conditions of the lemma.
\end{proof}

Lemma~\ref{thm:dressdel} shows that under suitable conditions,
a dressed Delaunay immersion is ambient isometric to the original Delaunay
immersion.

\begin{lemma}\label{thm:dressdel}
  Let $A\in\Loopsl{1}$ satisfy $A=A^\ast$,
  $\Phi=\exp(\zeta A)$ and $f_\lambda=\Sym{\lambda}{\Uni{r}{\Phi}}$
  the Delaunay associate family.
  Let $C\in\LoopSL{r}$, and suppose that
  $C$ is the boundary of an analytic map
  $C:\{r<\abs{\lambda}<1+\epsilon\}\to\mattwo(\bbC)$
  for some $\epsilon\in\bbR^{+}$ such that
  $\{\det C=0\}\subset\bbS^1$.
  Suppose that $C\exp(2\pi i A)C^{-1}$ satisfies the reality condition on
  $\bbS^1\setminus\{\det C=0\}$.
  Let $\widetilde{f}_\lambda=\Sym{\lambda}{\Uni{r}{C\Phi}}$.
  Then (i)
  There exists $\widetilde{A}$ of the form~\eqref{eq:del},
  $U_0\in\LoopuSL{r}$ and $B_0\in\LooppGL{r}$ such that
  $C\Phi = U_0 \exp(\zeta \widetilde{A} ) C_{+}$.
  (ii)
  Then there exists $c\in\bbR^{+}$
  and an isometry $T$ of $\matsu_2$ such that
  $\widetilde{f}_\lambda(\zeta)=T(f_\lambda(\zeta+c))$.
\end{lemma}

\begin{proof}
  Let $C_uC_{+}$ be the $r$-Iwasawa factorization of $C$.
  Because $C_u$ is analytic on $\calA_r$ with $\det C_u\in\matU_1$,
  $C_{+}$ is the boundary of an analytic map
  $C_{+}:\{r<\abs{\lambda}<1+\epsilon\}\to\mattwo(\bbC)$
  such that $\{\det C_{+}=0\}\subset\bbS^1$.
  By lemma~\ref{thm:unitcomm}, there exist
  analytic maps $U\in\LoopuSL{1}$ and $R:\bbS^1\to\mattwo(\bbC)$
  such that $C_{+}=UR$ and $[R,\,A]=0$.
  $U$ and $R$ can be extended to $\calA_{s}$ for some $s\in(r,\,1)$.

  Then $\widetilde{A}=UAU^{-1}=C_{+}AC_{+}^{-1}$ on $\{s<\abs{\lambda}<1\}$.
  But $UAU^{-1}$ extends analytically to $\calA_{s}$, and
  $C_{+}AC_{+}^{-1}$ extends holomorphically to $\{0<\abs{\lambda}<1\}$
  and meromorphically to $0$ with a simple pole in the upper right entry.
  Moreover, this extension satisfies $\widetilde{A}=\widetilde{A}^\ast$,
  since $UAU^{-1} = (UAU^{-1})^\ast$.
  It follows that $\widetilde{A}$ is of the form~\eqref{eq:del}.

  On $\calC_{r}$,
  \[
  C\Phi = C_u \exp(\zeta \widetilde{A} ) C_{+}.
  \]
  Hence $\Sym{\lambda}{\Uni{r}{C\Phi}}$ and
  $\Sym{\lambda}{\Uni{r}{\exp(\zeta\widetilde{A})}}$ are the
  same surface up to rigid motion. The result follows by
  lemma~\eqref{thm:del}(iv).
\end{proof}

\subsection{Delaunay asymptotics}

The following lemma estimates the growth rate of the gauge $B$ which gauges
the Maurer-Cartan form for the Delaunay associate family to the Delaunay
DPW potential.
Since the estimate is for $\abs{\lambda}$ near $1$, the explicit
Delaunay frame is not required; the growth rate can be estimated
by using the periodicity in the axial direction.
The result is essentially that $B(x+iy)$ grows in the axial direction
$x$ like $e^{c\abs{x}}$, where $c$ is the maximum value of the
Delaunay eigenvalue on $\bbS^1$.
The estimate is used in the asymptotics theorem~\ref{thm:asymptotic}
showing that a perturbation of the DPW Delaunay potential is asymptotically
Delaunay.

In the following, $\abs{X}$ denotes the matrix $2$-norm,
and for $r\in(0,\,1]$,
\[
\norm{X(\lambda)}_r = \max_{\lambda\in\bbC_r}\abs{X(\lambda)}.
\]

\begin{lemma}
  \label{thm:delasm}
  Let $\Sigma=\bbC$. Let $A$ be a Delaunay residue (equation~\eqref{eq:del}),
  let $\Phi=\exp(\zeta A)$, let $C\in\LoopprX$ and let
  $C\Phi=FB$ be the $r_0$-Iwasawa
  factorization of $C\Phi$ for some $r_0\in(0,1]$,
  and extend $B$ to $\Sigma\times\calD_1$ as in lemma~\ref{thm:rindependent}.
  Let $\mu$ be an eigenvalue of $A$ and let $c=\norm{\re\mu}_{1}$.
  Then there exists $c_0\in\bbR^{+}$ such that for all $\epsilon>0$,
  there exists $R(\epsilon)\in(0,\,1)$ such that
  \begin{equation*}
    \norm{B(\zeta,\,\lambda)}_r \le c_0\exp( (c+\epsilon)\abs{\re \zeta} )
  \end{equation*}
  for all $\zeta\in\bbC$ and all $r\in(R(\epsilon),\,1]$.
\end{lemma}

\begin{proof}
  First we prove the theorem in the case $C=\id$.
  With $\zeta=x+iy$,
  the screw symmetry of the Delaunay family implies that
  $F$ decouples into $x$- and $y$-dependent factors
  $F=\exp(i y A)F_1(x)$ for some $F_1:\Sigma\to\LoopuSL{1}$.
  Then $B(x)=F_1(x)^{-1}\exp(x A)$ can be estimated by
  estimating $\exp(x A)$ and $F_1(x)$.

  \emph{Step 1: estimate $\exp(xA)$.}
  From the formula
  \[
  \exp(xA)=\half e^{x\mu}(\id+\mu^{-1}A)+
  \half e^{-x\mu}(\id-\mu^{-1}A)
  \]
  we obtain the pointwise estimate
  \[
  \abs{\exp(xA)} \le
  (\max\abs{\id\pm\mu^{-1}A}) \exp(\abs{\re\mu}\abs{x}),
  \]
  for all $x\in\bbR$ and all $\lambda$ at which
  $\max\abs{\id\pm\mu^{-1}A}$ is finite.
  Since $\max\norm{\id\pm \mu^{-1}A}_r$ is continuous and finite
  at $r=1$, there exists $R_1\in(0,\,1)$ such that
  $\max\norm{\id\pm \mu^{-1}A}_r$ is finite for all
  $r\in[R_1,\,1]$. Then
  \[
  c_1=\sup_{r\in[R_1,\,1]}\max(\norm{\id\pm\mu^{-1}A}_r).
  \]
  is finite.
  Then for all $x\in\bbR$ and all $r\in(R_1,\,1]$,
  \begin{equation*}
    \norm{\exp(xA)}_r \le c_1 \exp(\norm{\re\mu}_{r}\abs{x}).
  \end{equation*}

  The continuity of $\norm{\re\mu}_{r}$ at $r=1$
  together with $\norm{\re\mu}_{1}=c$
  imply that all $\epsilon>0$ there exists
  $R\in(0,\,1)$ such that
  for all $r\in(R,\,1]$,
  $\norm{\re\mu}_{r} < c+\epsilon$.
  Hence for all $\epsilon>0$ there exists $R\in(0,\,1)$ such that
  for all $x\in\bbR$ and all $r\in(R,\,1]$,
  \begin{equation}
    \label{eq:normExp}
    \norm{\exp(xA)}_r \le c_1 \exp((c+\epsilon)\abs{x}).
  \end{equation}

  \emph{Step 2: estimate $F_1(x)$.}
  $F_1$ is periodic in axial direction the sense that
  there exist $\rho\in\bbR$ and $M\in\LoopuSL{1}$ such that
  \[
  F_1(x_0+n\rho) = M^n F(x_0)
  \]
  for all $x\in\bbR$ and all $n\in\bbZ$.
  There exists $R_2\in(0,\,1)$ such that
  $\norm{F(x,\,\lambda)}_r$ is finite for all $r\in[R_2,\,1]$
  and all $x\in\bbR$. Then
  \[
  c_2 = \sup_{(x,\,r)\in[0,\rho)\times[R_2,\,1]}
  \norm{F(x,\,\lambda)}_r
  \]
  is finite.

  Given any $x\in\bbR$, there exists $x_0\in[0,\,\rho)$ and
  $n\in\bbZ$ such that $x=x_0+n\rho$.
  Hence
  \[
  \norm{F_1(x)}_r \le c_2(\norm{M}_r)^n.
  \]
  The continuity of $\norm{M}_{r}$ at $r=1$ together with
  $\norm{M}_{1}=1$ imply that
  for all $\epsilon'>0$ there exists $R$ such that for all $r\in(R,\,1)$,
  $\norm{M}_{r}<1+\epsilon'$.
  Given $\epsilon>0$,
  let $\tilde{\epsilon}=\min(\epsilon,\,1)$, and
  choose $\epsilon'=\exp(\rho\tilde{\epsilon})$.
  Then there exists $R$ such that for all $r\in(R,\,1)$,
  $\norm{M}_{r}<1+\epsilon' = \exp(\rho\tilde{\epsilon})$.
  Hence
  \[
  (\norm{M}_{r})^\abs{n} < \exp(\rho\tilde\epsilon\abs{n})=
  \exp(\tilde\epsilon\abs{x-x_0})\le
  \exp(\tilde\epsilon\rho)\exp(\tilde\epsilon\abs{x}).
  \]
  Hence with $c_3 = c_2\exp(\rho)$,
  \begin{equation}
    \label{eq:normF}
    \norm{F_1(x)}_r \le c_3 \exp(\epsilon\abs{x}).
  \end{equation}

  \emph{Step 3: estimate $B$.} 
  $B(x)=F_1^{-1}(x)\exp(x A)$, so
  \[
  \norm{B(x,\,\lambda)}_r \le
  \norm{F_1(x,\,\lambda)}_r\norm{\exp(xA(\lambda))}_r.
  \]
  Given $\epsilon>0$, by~\eqref{eq:normExp} and
  and~\eqref{eq:normF} we can choose $R$ such that
  for all $x\in\bbR$ and all $r\in(R,\,1]$,
  \[
  \norm{\exp(xA)}_r \le c_1 \exp((c+\epsilon/2)\abs{x})
  \]
  and
  \[
  \norm{F_1(x)}_r \le c_3 \exp((\epsilon/2)\abs{x}).
  \]
  Then with  $c_4 = c_1 c_3$,
  \[
  \norm{B(x,\,\lambda)}_r \le c_4 \exp((c+\epsilon)\abs{x}).
  \]

  Now we prove the theorem for general $C$.
  By theorem~\ref{thm:dressdel}, $C\Phi=U_0 \widetilde{\Phi} B_0$,
  where $\widetilde{\Phi}=\exp(\zeta\widetilde{A}$,
  $U_0\in\LoopuSL{r}$, $B_0\in\LoopGL{r}$.
  Let $\widetilde{\Phi}=\widetilde{F}\widetilde{B}$ be the Iwasawa
  factorization of $\widetilde{\Phi}$.  Then
  $C\Phi=(U_0\widetilde{F})(\widetilde{B}B_0)$ is the Iwasawa
  factorization of $C\Phi$.

  Then since $B_0$ is $\zeta$-independent,
  for any $\epsilon>0$, there exists $R\in(0,\,1)$ such that
  for all $x\in\bbR$ and all $r\in(R,\,1]$,
  \[
  \norm{\widetilde{B}(x,\,\lambda)B_0(\lambda)}_r \le 
  \norm{\widetilde{B}(x,\,\lambda)}_r \norm{B_0(\lambda)}_r \le
  c_0 \exp((c+\epsilon)\abs{x}),
  \]
  where $c_0 = c_4 \sup_r\norm{B_0}_r$.
\end{proof}

\begin{figure}[ht]
  \centering
  \includegraphics[width=175pt]{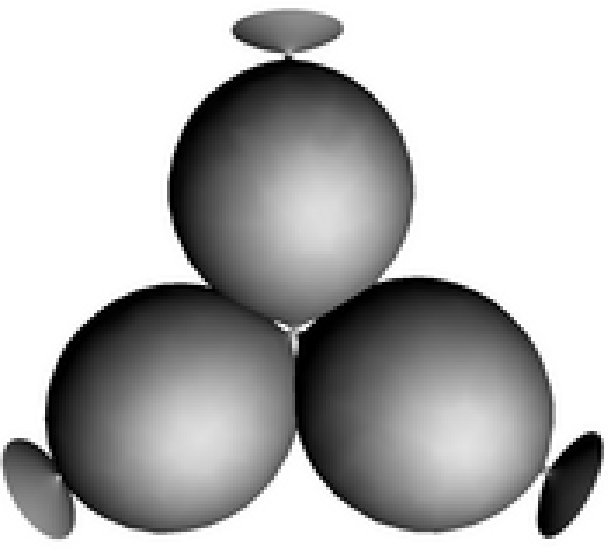}
  \includegraphics[width=175pt]{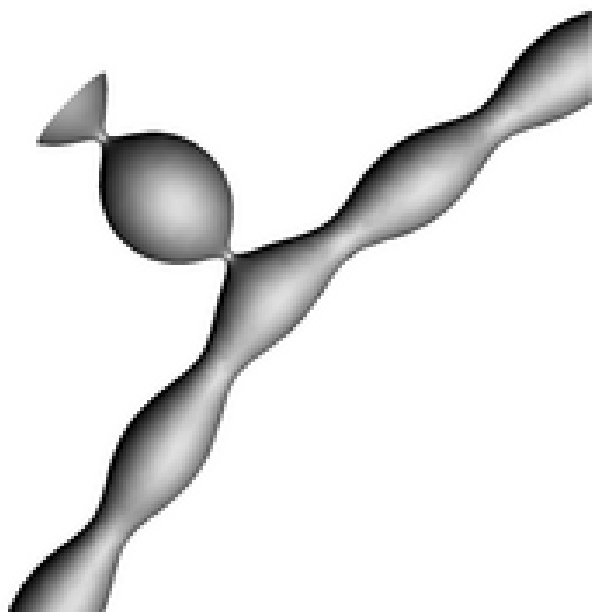}
\caption{Trinoids with small necks.
The lobes of the left example (necksizes 
$\bigl(\tfrac{1}{40},\,\tfrac{1}{40},\,\tfrac{1}{40}$)
intersect, making the surface Alexandrov embedded.
The example on the right (necksizes
$\bigl(\tfrac{1}{3},\,\tfrac{1}{3},\,\tfrac{1}{12}\bigr)$
can be viewed as a Delaunay surface with a small-necked Delaunay
end added. The Delaunay surface bends slightly to balance the added end.}
\end{figure}

\typeout{== perturb.tex =============================================}
\section{Perturbations of Delaunay immersions}
\label{sec:asymptotics}

\subsection{Perturbations at a simple pole}
\label{sec:zap}

The following lemma extends a basic result in ODE theory to the
context of loops.

\begin{lemma}
  \label{thm:zap}
  Let $r\in(0,\,1)$.
  Let $\xi_0,\,\xi\in\Potentialsl{\Sigma^\ast}{r}$ be potentials with
  expansions in $z$ at $z=0$
  \[
  \xi_0 = A\frac{dz}{z},\quad \xi = A\frac{dz}{z} + Bdz + \Order(z^1)dz.
  \]
  Let $\mu$ be an eigenvalue of $A$ and suppose that either
  \begin{enumerate}
  \item $\mu\notin\half\bbZ^\ast$ along $\bbC_r$, or
  \item $\mu\notin\half\bbZ^\ast\setminus\{\pm\half\}$ along $\bbC_r$
    and $[A,\,B]=0$.
  \end{enumerate}
  Then in a neighborhood $U$ of $p$ there exists a unique analytic map
  $P:U\times\to\LoopSL{r}$
  such that
  \begin{equation}
    \label{eq:P}
    \xi = \gauge{\xi_0}{P},\quad P(0,\,\lambda)=\id.
  \end{equation}
  In the case $[A,\,B]\equiv 0$, $P=\id+Bz+\Order(z^2)$.
\end{lemma}

\begin{proof}
  In case (i), a unique solution to~\eqref{eq:P}
  exists by a standard result in ODE theory,

  In the case (ii), if $\mu(\lambda_0)=\pm\half$,
  a calculation of the series $P=\sum_{k=0}^\infty P_k z^k$ shows that
  the $P_k$ are holomorphic in $\calC_r$, and $P_1=B$.
\end{proof}

\subsection{Gauging away the constant term}

Lemma~\ref{thm:gaugeconst} constructs a gauge and coordinate change
which removes the constant term in a perturbed Delaunay potential.

\begin{lemma}\label{thm:gaugeconst}
  Let $r\in(0,\,1]$.
  Let $\Sigma$ be a Riemann surface and $p\in\Sigma$.
  Let $\xi\in\Potentialsl{\Sigma}{r}$ with expansion
  \begin{equation}
    \label{eq:xig}
    \gauge{\xi}{g} =
    \xi_{-1}\frac{dz}{z} + \xi_0 dz + \Order(z^1)\,dz.
  \end{equation}
  Let $\mu$ be an eigenvalue of $\xi_{-1}$ and suppose
  \begin{enumerate}
  \item
    $\res_{\lambda=0}\mu^2 \ne 0$
  \item
    for every $\lambda_0\in\calD_{r}$, if $\mu(\lambda_0)\in\{\pm\half\}$,
    then $\xi_0|_{\lambda=\lambda_0}=0$.
  \end{enumerate}
  Then there exists a neighborhood $U\in\Sigma$ of $p$,
  an analytic map $g:U\cross\calD_1\to\matGL_2(\bbC)$
  such that $g(z,\,0)$ takes values in $\TriGL$, and a
  conformal coordinate $\tilde{z}:U\to\bbC$ with $\tilde{z}(p)=0$
  such that the expansion of $\gauge{\xi}{g}$ in $\tilde{z}$ at
  $\tilde{z}=0$ is
  \begin{equation}
    \label{eq:xig2}
    \gauge{\xi}{g} =
    \xi_{-1}\frac{d\tilde{z}}{\tilde{z}} +
    \Order(\tilde{z}^1)\,d\tilde{z}.
  \end{equation}
\end{lemma}

\begin{proof}
  For any $k\in\bbC$, define $g_1$ and $g$ by
  \begin{equation*}
    \begin{split}
      &u = 4\mu^2-1,\quad
      v\id = \xi_{-1}\xi_{0}+\xi_{0}\xi_{-1},\\
      &g_1 = (k-2u^{-1}v)\xi_{-1} +
      u^{-1}(\xi_0-[\xi_{-1},\,\xi_{0}]),\quad
      g = \id + g_1 z.
    \end{split}
  \end{equation*}
  A calculation shows that
  \begin{equation*}
    (\id+\ad_{\xi_{-1}})g_1 = k\xi_{-1} -\xi_0,
  \end{equation*}
  from which it follows that
  \begin{equation*}
    \gauge{\xi}{g} = \xi_{-1}\frac{dz}{z} + k\xi_{-1}\,dz + \Order(z^1)dz.
  \end{equation*}

  Assumption (i) implies that $u^{-1}v$ is holomorphic at
  $\lambda=0$ so $k = \lim_{\lambda\to 0}2u^{-1}v$ exists and is finite.
  With this choice of $k$, a calculation shows $g_1$ is holomorphic at
  $\lambda=0$.
  Assumption (ii) implies that $g_1$ is holomorphic in $\calD_r^\ast$,
  and hence in $\calD_r$.
  A calculation shows that $g_1(0)\in\TriGL$ and hence $g(z,\,0)$
  takes values in $\TriGL$.

  Since $g(0,\,\lambda)=\id$, a continuity argument shows that
  $\det{g}\ne 0$ in a sufficiently small neighborhood of $z=0$.
  In the coordinate $\tilde{z}$ defined by
  $z = \tilde{z} - k \tilde{z}^2$ in a neighborhood of $z=0$,
  $\xi$ has the expansion~\eqref{eq:xig2}
\end{proof}

\subsection{Monodromy at simple poles}

The following lemma computes the eigenvalues of the monodromy of a
perturbed potential $\xi$ at a simple pole in terms of the
residue of $\xi$.

\begin{lemma}\label{thm:eigen}
  Let $r\in(0,\,1)$.
  Let $\xi\in\Potentialsl{\Sigma^\ast}{r}$ be a potential with
  expansions in $z$ at $z=0$
  $\xi = Adz/z + \Order(z^0)dz$,
  and let $\mu$ be an eigenvalue of $A$.
  Suppose $\xi$ satisfies condition (i) or (ii) of lemma~\ref{thm:zap}.
  Then

  (i) 
  If $\Phi:\widetilde{\Sigma}\times\calC_r\to\matGL_2(\bbC)$
  is a solution to the ODE $d\Phi=\Phi\xi$,
  and $M$ is the monodromy of $\Phi$ at $z=0$,
  then the eigenvalues of $M$ are $\exp(\pm 2\pi i \mu)$.

  (ii) 
  If $\Phi:\widetilde{\Sigma}\times\calC_r\to\mattwo(\bbC)$
  is a solution to the ODE $d\Phi=\Phi\xi$ with
  $\det\Phi\notequiv 0$, and
  the monodromy $M$ of $\Phi$ at $z=0$
  extends analytically to $\calC_r$
  across $\{\det\Phi=0\}$, then the eigenvalues of $M$ are
  $\exp(\pm 2\pi i \mu)$.
\end{lemma}

\begin{proof}
  Since (ii) implies (i) we prove (ii).
  Let $\xi_0 = A dz/z$.
  By lemma~\ref{thm:zap} there exists a unique analytic map
  $P:U\to\LoopSL{r}$
  such that $\xi=\gauge{\xi_0}{P}$ and $P(0,\,\lambda)=\id$.
  Then there exists an analytic map $C:\calC_r\to\mattwo(\bbC)$ such that
  $\Phi = C\exp((\log(z) A)P$.

  Then $M = C\exp(2\pi i A)C^{-1}$ on $\calC_r\setminus\{\det\Phi = 0\}$,
  so the eigenvalues of $M$ are $\exp(2\pi i \mu)$ 
  on $\calC_r\setminus\{\det\Phi = 0\}$.
  Since by hypothesis $M$ extends analytically to $\calC_r$,
  the eigenvalues of $M$ extend analytically to $\calA_r$,
  and hence are $\exp(2\pi i \mu)$ on $\calC_r$.
\end{proof}

\subsection{Perturbed Delaunay asymptotics}
\label{sec:del-pert}

In this section it is shown that the immersion obtained from a suitable
perturbation of a Delaunay potential is asymptotic to the base half-Delaunay
surface (theorem~\ref{thm:asymptotic}).

In the following, $\abs{X}$ denotes the matrix $2$-norm,
and for a compact set $S\subset\bbC^\ast$,
\[
\norm{X(\lambda)}_S = \max_{\lambda\in S}\abs{X(\lambda)}.
\]

The asymptotics theorem below shows that under certain conditions,
the CMC immersion produced by a perturbation of a Delaunay potential is
asymptotic to a half Delaunay surface.

\begin{theorem}[Delaunay asymptotics theorem]
  \label{thm:asymptotic}
  Let $\Sigma$ be a punctured annular neighborhood of $0$ and
  $\widetilde{\Sigma^\ast}\to\Sigma^\ast$ its universal cover.
  Let $\xi_\Del=A\,dz/z\in\Potentialsl{\Sigma^\ast}{r}$ where
  $A\in\Loopsl{r}$ is of the form~\eqref{eq:del}.
  Let $\mu$ be an eigenvalue of $A$, let $k\in\bbZ$, $k\ge 1$, and suppose
  $\norm{\re\mu}_{\bbS^1}<k$.
  Let $\xi\in\Potentialsl{\Sigma^\ast}{r}$ be a perturbation of $\xi_\Del$
  whose expansion of $\xi$ at $z=0$ is
  \[
  \xi=Az^{-1}\,dz + \Order(z^{2k-1})dz.
  \]
  Let $\Phi:\widetilde{\Sigma^\ast}\times\calA_r\to\mattwo(\bbC)$,
  with values in $\matGL_2(\bbC)$ off $\bbS^1$,
  satisfy $d\Phi=\Phi\xi$.
  Let $\Phi_\Del =\Phi P^{-1}$, where $P$ is the gauge of lemma~\ref{thm:zap}
  with $\gauge{\xi_\Del}{P}=\xi$ and $P(0,\,\lambda)=\id$.
  Let $f_\Del=\Sym{}{\Uni{s}{\Phi_\Del}}$.
  Then
  \begin{align}
    \label{eq:asymf}
    &\lim_{z\rightarrow 0}\norm{f-f_\Del}_{\bbS^1}=0\\
    \label{eq:asymdf}
    &\lim_{z\rightarrow 0}\norm{df-df_\Del}_{\bbS^1}=0.
  \end{align}
\end{theorem}

\begin{proof}
  Let $\Phi_\Del = F_\Del B_\Del$ and $\Phi=FB$ be the $r$-Iwasawa factorizations
  of $\Phi_\Del$ and $\Phi$ respectively. By hypothesis, the monodromy
  of $\Phi_\Del$ is $r$-unitary. It follows that the monodromy of $\Phi$ is
  $r$-unitary, and that $B_\Del$, $B$ and $F_\Del^{-1}F$ are monodromy-free on
  $\Sigma$.

  By lemma~\ref{thm:zap} the expansion of $P(z)$ at $z=0$ is
  \[
  P(z) = \id + \msmall{\sum_{j=2k}^{\infty}}P_j z^j.
  \]
  Then
  \begin{equation*}
    B_\Del PB_\Del ^{-1}-\id = \msmall{\sum_{j=2k}^{\infty}}B_\Del P_j B_\Del^{-1} z^j,
  \end{equation*}
  so
  \begin{equation*}
    \norm{B_\Del PB_\Del ^{-1}-\id}_{\calC_r} \le
    \msmall{\sum_{j=2k}^{\infty}}\norm{B_\Del}_{\calC_r}\norm{P_j}_{\calC_r}
    \norm{B_\Del^{-1}}_{\calC_r}\abs{z}^j.
  \end{equation*}
  By hypothesis $c=\norm{\re\mu}_{\bbS_1}<k$.
  Let $\epsilon\in(0,\,k-c)$, so that
  $l = c+\epsilon < k$. By lemma~\ref{thm:delasm},
  there exists $R$ such that
  for all $r\in(R,\,1]$ and all $z\in\Sigma$,
  $\abs{B_\Del} < c_0 \abs{z}^{-l}$
  for some constant $c_0\in\bbR^{+}$.
  Then
  \begin{equation*}
    \norm{B_\Del PB_\Del ^{-1}-\id}_{\calC_r} \le
    \msmall{\sum_{j=2k}^{\infty}}\norm{P_j}_{\calC_r}\abs{z}^{j-2l}.
  \end{equation*}
  By the choice of $l$, the exponent $j-2l>0$ for all $j\ge 2k$, so
  \begin{equation*}
    \lim_{z\rightarrow 0}\norm{B_\Del PB_\Del^{-1}-\id}_{\calC_r}=0.
  \end{equation*}
  The holomorphicity of $B_\Del PB_\Del ^{-1}$ in $\lambda$ with Cauchy's integral
  formula implies
  \begin{equation*}
    \lim_{z\rightarrow 0}\norm{ (B_\Del PB_\Del^{-1})'}_{\calC_r}=0.
  \end{equation*}
  With $G = F_\Del^{-1}F = \Uni{r}{B_\Del PB_\Del ^{-1}}$,
  \begin{align}
    \label{eq:asymG}
      &\lim_{z\rightarrow 0}\norm{G-\id}_{A} = 0\\
    \label{eq:asymGp}
      &\lim_{z\rightarrow 0}\norm{G'}_{A} = 0\\
    \label{eq:asymB}
      &\lim_{z\rightarrow 0}\norm{BB_\Del^{-1}-\id}_{D} = 0
  \end{align}
  for every compact subset $A\subset\calA_r$ and $D\subset\calD_r$.

  From the Sym formula~\eqref{eq:sym} we get
  \begin{equation*}
    f-f_\Del = -2\abs{H}^{-1}F_\Del G' F^{-1}.
  \end{equation*}
  Then since
  $\norm{F_\Del}_{\bbS^1}=1$ and $\norm{F^{-1}}_{\bbS^1}=1$
  we have
  \begin{equation*}
    \norm{f-f_\Del}_{\bbS^1} \le 2\abs{H}^{-1}\norm{G'}_{\bbS^1},
  \end{equation*}
  and equation~\eqref{eq:asymf} follows.

  \Note{Check.}
  Differentiating the Sym formula~\eqref{eq:sym} yields
  \[
  df = -2H^{-1}F\Theta'F^{-1},\quad
  df_\Del=-2H^{-1}F_\Del\Theta_\Del'F_\Del^{-1},
  \]
  where
  \begin{equation*} 
    \Theta' = \fourth H v E,\quad
    \Theta_\Del' = \fourth H v_\Del E_\Del,
  \end{equation*}
  $v^2 dz\tensor d\ol{z}$ and $v_\Del^2 dz\tensor d\ol{z}$
  are the metrics of $f$ and $f_\Del$ respectively,
  \begin{equation*}
    E= -i{\abs{\alpha}}^{-1}
    \begin{pmatrix}
      0 & \alpha\lambda^{-1} \\ \ol{\alpha}\lambda & 0
    \end{pmatrix},\quad
    E_\Del= -i{\abs{\alpha_\Del}}^{-1}
    \begin{pmatrix}
      0 & \alpha_\Del\lambda^{-1} \\ \ol{\alpha_\Del}\lambda & 0
    \end{pmatrix},
  \end{equation*}
  and $\alpha$, $\alpha_\Del$ are defined by
  \[
  \xi = 
  \begin{pmatrix}
    0 & \alpha \\ 0 & 0
  \end{pmatrix}\lambda^{-1} + \Order(\lambda^0)dz,\quad
  \xi_\Del = 
  \begin{pmatrix}
    0 & \alpha_\Del \\ 0 & 0
  \end{pmatrix}\lambda^{-1} + \Order(\lambda^0)dz.
  \]
  Then
  \begin{equation*}
    df-df_\Del = -\half v_\Del F_\Del(v_\Del^{-1}v G E - E_\Del G)F^{-1}.
  \end{equation*}
  Then since
  $\norm{F_\Del}_{\bbS^1}=1$ and $\norm{F^{-1}}_{\bbS^1}=1$
  we have
  \begin{equation*}
    \begin{split}
      &\norm{df-df_\Del}_{\bbS^1} \\
      &\quad\quad\le \half \norm{v_\Del}_{\bbS^1}
      \norm{v_\Del^{-1}v G E - E_\Del G}_{\bbS^1}\\
      &\quad\quad\le \half \norm{v_\Del}_{\bbS^1} \left(
        \norm{v_\Del v^{-1}G-\id}_{\bbS^1}\norm{E}_{\bbS^1}
        + \norm{G-\id}_{\bbS^1}\norm{E_\Del}_{\bbS^1} +
        \norm{E-E_\Del}_{\bbS^1}\right).
    \end{split}
  \end{equation*}
  $\norm{v_\Del}_{\bbS^1}$,
  $\norm{E}_{\bbS^1}$ and $\norm{E_\Del}_{\bbS^1}$ are bounded on $\Sigma$,
  and
  \[
  \lim_{z\to 0}\norm{E-E_\Del}_{\bbS^1}= 0.
  \]
  By equation~\eqref{eq:asymB} and remark~\ref{rem:dpw},
  \begin{equation*}
      \lim_{z\to 0}v_\Del^{-1}v=1,
  \end{equation*}
  from which it follows, using equation~\eqref{eq:asymG}, that
  \[\lim_{z\to 0}\norm{v_\Del v^{-1}G-\id}_{\bbS^1}= 0.\]
  Equation~\eqref{eq:asymdf} follows.
\end{proof}

\begin{corollary}
  If in theorem~\ref{thm:asymptotic} $A$ satisfies
  equations~\eqref{eq:del}--\eqref{eq:del-close}, the weight $w$
  associated to $A$ satisfies $w>-3$, the expansion of $\xi$ is
  \[
  \xi=Az^{-1}\,dz + \Order(z^1)dz,
  \]
  and the monodromy $\Phi$ at $z=0$ is in $\LoopuSL{r}$, by
  lemma~\ref{thm:dressdel} $f_\Del$ is a Delaunay associate family with
  weight $w$, the Delaunay and perturbed surfaces are closed at
  $\lambda=1$, and the theorem shows $C^1$ convergence of the
  perturbed surface to the Delaunay surface.
\end{corollary}

\begin{figure}[ht]
  \centering
  \includegraphics[width=175pt]{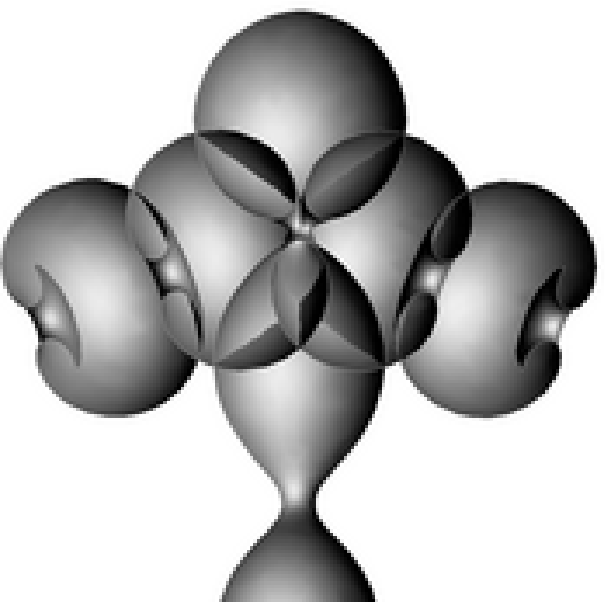}
  \includegraphics[width=175pt]{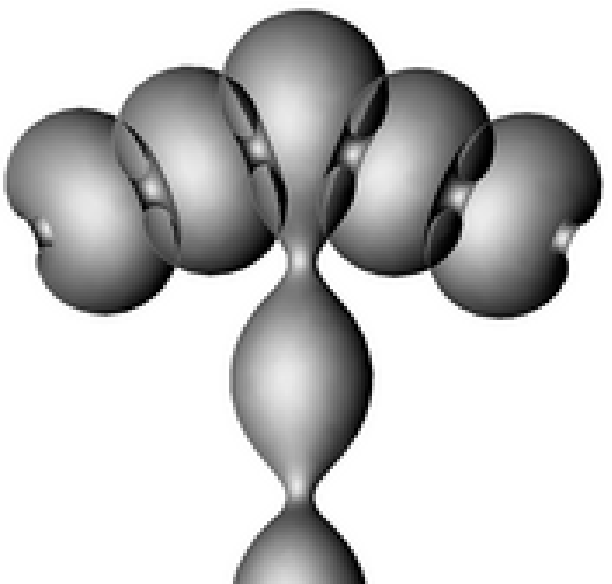}
\caption{A pair of CMC trinoids with two nodoid ends (necksizes
$\bigl(\tfrac{1}{6},\,-\tfrac{1}{6},\,-\tfrac{1}{6}\bigr)$).}
\end{figure}

\typeout{== unitarize.tex =============================================}
\section{Unitarization of loop group monodromy representations}
\label{sec:glue}

This section proves the ``gluing theorem'' (theorem~\ref{thm:glue}):
if a monodromy representation of the ODE~\eqref{eq:ode} on the $n$-punctured
Riemann sphere is unitarizable
pointwise on $\bbS^1$, then the monodromy representation is
unitarizable by a dressing matrix on an $r$-circle which is
analytic in $\lambda$.  The proof is based on
lemmas~\ref{thm:birkhoff1}--\ref{thm:section}.

\begin{notation}
  \label{not:unit}
  $M\in\matGL_2(\bbC)$ is \emph{unitarizable} if there exists
  $C\in\matGL_2(\bbC)$ such that $CMC^{-1}\in\matU_2$.

  The set $\calM=\{M_1,\dots,M_n\}\subset\matGL_2(\bbC)$
  is \emph{simultaneously unitarizable} iff
  for all $j\in\{1,\dots,n\}$ there exists
  $C\in\matGL_2(\bbC)$ such that
  $CM_jC^{-1}\in\matU_2$.

  $\calM$ is \emph{nondegenerate} iff $[M_i,\,M_j]\ne 0$ for some pair
  $i\ne j$.
\end{notation}

\subsection{Birkhoff factorizations}
\label{sec:birkhoff}

Two Birkhoff factorizations are given for singular loops on
$\bbS^1$: a scalar version (lemma~\ref{thm:birkhoff1}) and a matrix version
(lemma~\ref{thm:birkhoff2}).
\medskip

\begin{lemma}\label{thm:birkhoff1}
  Let $f:\bbS^1\to\bbR^{\ge 0}$ be an analytic map 
  with $f\notequiv 0$.
  Then there exists an analytic map
  $h:\bbS^1\to\bbC$ which is the boundary of an
  analytic map $\calD_1\to\bbC^\ast$, such that $f=h^\ast h$.
\end{lemma}

\begin{proof}
  Since $f$ is real and non-negative, each of its zeros is of even order.
  Let $\{a_1,\dots,a_n\}\subset\bbS^1$ be the zeros of $f$, each with
  multiplicity two, and let $q=\prod_{j=1}^n(\lambda-a_j)$.
  Then the function $g = f/(q^\ast q)$ has no zeros on $\bbS^1$ and
  satisfies $g=g^\ast$. Let
  \[
  g=r\lambda^p g_{-}g_{+}
  \]
  be the (rank 1) Birkhoff factorization of
  $g$, such that $g_{+}$ extends analytically without zeros to
  $\overline\calD_1$,
  $g_{-}$ extends analytically without zeros to $\overline\calD_1'$,
  and normalized with $r\in\bbC$, $g_{+}(0)=1$ and $g_{-}(\infty)=1$.
  But $g^\ast=g$ on $\bbS^1$, so on $\bbS^1$ we have the equality
  \begin{equation*}
    r\lambda^p g_{-}g_{+} = \overline{r}\lambda^{-p}g_{+}^\ast g_{-}^\ast.
  \end{equation*}
  By the uniqueness of the Birkhoff factorization,
  $g_{-}=g_{+}^\ast$, $p=0$ and $r=\overline{r}$. Since $f$ is nonnegative
  on $\bbS^1$, $r$ is positive. Then the function
  \[
  h=\sqrt{r}g_{+}q
  \]
  is analytic on $\bbS^1$, is the boundary of the map $h:\calD_1\to\bbC^\ast$
  and satisfies $f = h^\ast h$.
\end{proof}

\begin{lemma}\label{thm:birkhoff2}
  Let
  $X:\bbS^1\to\mattwo(\bbC)$ be a
  positive semidefinite analytic map
  with $\det{X}\notequiv 0$.
  Then there exists
  $C\in\LoopprX$ and an analytic map $f:\bbS^1\to\bbR^{\ge 0}$
  such that $fX=\left.C^\ast C\right|_{\bbS^1}$.
\end{lemma}

\begin{proof}
  The map $X$ can be written
  \begin{equation*}
    X=
    \begin{pmatrix}
      x_1 & y \\ y^\ast & x_2
    \end{pmatrix}
  \end{equation*}
  where the functions $x_1,\,x_2$ satisfy $x_1=x_1^*$ and $x_2=x_2^*$,
  are real-valued and non-negative on $\bbS^1$, and $x_1\notequiv 0$,
  $x_2\notequiv 0$ on $\calA_r$.

  The function $d=\det X$ satisfies $\det X\notequiv 0$
  on $\calA_r$, and since $X$ is positive semidefinite, $d$ is
  real-valued and non-negative on $\bbS^1$.
  $d=e^\ast e$ be the singular Birkhoff factorizations of $d$
  (lemma~\ref{thm:birkhoff1}).  Let
  \begin{equation*}
    Y=
    \begin{pmatrix}
      x_1 & y \\ 0 & e
    \end{pmatrix}.
  \end{equation*}
  Then $Y$ is a analytic map on $\bbS^1$ which satisfies
  \[
  x_1 X=Y^\ast Y.
  \]

  For some $r\in(0,\,1)$, $X$ extends analytically to
  a map $\widetilde{X}:\calA_r\to\mattwo(\bbC)$
  such that $\widetilde{X}_{11}$ and $\det\widehat{X}$ have no zeros in
  $\calA_s\setminus\bbS^1$.
  Then $Y$ likewise extends analytically to a map
  $\widetilde{Y}:\calA_r\to\mattwo(\bbC)$
  such that $\det\widehat{Y}$ have no zeros in
  $\calA_r\setminus\bbS^1$.
  Let $\widetilde{Y}|_{\calC_s}=Y_{u}Y_{+}$ be the $s$-Iwasawa factorization
  of $\widetilde{Y}|_{\calC_s}$ for any $s\in(r,\,1)$.
  Since $\widetilde{Y}|_{\calC_s}$ and $Y_u$ are the boundaries
  of analytic maps on $\calA_s$ with nonzero determinants on
  $\calA_s\setminus\bbS^1$,
  then $Y_{+}$ is the boundary of an analytic map
  $\widetilde{Y}_{+}:\calD_1\to\matGL_2(\bbC)$.
  Then $x_1 X = \widetilde{Y}_{+}^\ast \widetilde{Y}_{+}|_{\bbS^1}$, so
  $C=Y_{+}$ and $f=x_1$ are the required maps.
\end{proof}

\subsection{Holomorphic vector bundles and unitarization}
\label{sec:glue0}

We prove several pointwise and holomorphic lemmas relating to simultaneous
unitarization.

\begin{lemma}\label{thm:mindim}
  Let
  \begin{equation*}
    L_\lambda:\bbC^m\to\bbC^n
  \end{equation*}
  be a family of linear maps which depends analytically on
  $\lambda\in\bbC^\ast$. Let
  \begin{equation*}
    r=\min_{\lambda\in\bbC^\ast}\dim\ker L_\lambda.
  \end{equation*}
  Then (i) $\dim\ker L_\lambda=r$ on $\bbC^\ast\setminus P$ for some subset
  $P\in\bbC^\ast$ of isolated points, and
  (ii) there exists a trivial analytic rank-$r$ bundle
  $E\to\bbC^\ast$ such that
  $E_\lambda\subseteq\ker L_\lambda$ on $\bbC^\ast$, and
  $E_\lambda=\ker L_\lambda$ on $\bbC^\ast\setminus P$.
\end{lemma}

\begin{lemma}\label{thm:isu}
  Let $U_1,\,U_2\in\matU_2$ with $[U_1,\,U_2]\ne 0$.
  Let $A\in\matGL_2(\bbC)$ and suppose that
  $AU_1A^{-1}\in\matU_2$ and
  $AU_2A^{-1}\in\matU_2$.
  Then $A\in\bbR^+\times\matU_2$.
\end{lemma}
\begin{proof}
  Choose a basis for which $U_1$ is diagonal.
  Factor $A=UT$, where $U\in\matU_2$ and $T\in\TrirGL$.
  Then $TU_1T^{-1}\in\matU_2$ implies $T$ is diagonal,
  and $TU_2T^{-1}\in\matU_2$ implies $T\in\bbR^{+}\id$.
  Hence $A = UT\in\bbR^+\times\matU_2$.
\end{proof}

\begin{lemma}\label{thm:dim1}
  (1)
  Let $M_1\in\matGL_2(\bbC)\setminus\{\pm\id\}$ be unitarizable.
  Let
  $L_1:\mattwo(\bbC)\to \mattwo(\bbC)$ be the linear map
  defined by
  \begin{equation*}
    L_1(X) = XM_1-{M_1^\ast}^{-1}X.
  \end{equation*}
  Then $\dim\ker L_1=2$.

  (2)
  Let $M_1,\dots,M_n\in\matGL_2(\bbC)$, $n\ge 2$, and suppose that
  $\{M_1,\dots,M_m\}$
  is simultaneously unitarizable and nondegenerate.  Let
  $L:\mattwo(\bbC)\to (\mattwo(\bbC))^n$ be the linear map
  defined by
  \begin{equation*}
    L(X)=(XM_1-{M_1^\ast}^{-1}X,\dots,XM_n-{M_n^\ast}^{-1}X).
  \end{equation*}
  Then $\dim\ker L=1$.
\end{lemma}

\begin{proof}
  To show (1),
  by hypothesis there exists $C\in\matGL_2(\bbC)$ such that
  $CM_1C^{-1}\in\matSU_2$. Let $X_0=C^\ast C$.
  A calculation shows that $X_0\in\ker L_1$ iff $[X_0^{-1}X,\,M_1]=0$.
  Since the space of commutators with $M$ is $2$-dimensional,
  then $\dim\ker L_1=2$ and $\ker L_1 = \Span\{X_0,\,X_0 M_1\}$.

  To show (2),
  assume without loss of generality that $M_1\notin\{\pm\id\}$.
  By hypothesis there exists $C\in\matGL_2(\bbC)$ such that
  $CM_jC^{-1}\in\matSU_2$. Let $X_0=C^\ast C$.
  Then $X_0\in\ker L$ so $\dim\ker L\ge 1$.
  But $\ker L\subset\ker L_1$, so $\dim\ker L\le 2$.

  Suppose $\dim\ker L=2$. Then as above,
  $\ker L = \Span\{X_0,\,X_0 M_j\}$ for each $j$.
  Hence for all $i,\,j$,
  $X_0M_i\in\Span\{X_0,\,X_0 M_j\}$,
  so $M_i\in\Span\{\id,\,M_j\}$
  so $[M_i,\,M_j]=0$,
  contrary to the hypothesis of the lemma.
\end{proof}

\begin{notation}
  Let $E\to\bbS^1$ be a vector bundle.
  $E(\lambda)$ denotes the fiber of $E$ over $\lambda\in\bbS^1$.
  $E^\ast$ denotes the vector bundle whose fiber over $\lambda\in\bbS^1$ is
  \begin{equation*}
    \{\transpose{\overline{X}}\suchthat X\in
    E({{\overline{\lambda}}^{-1}})\}.
  \end{equation*}
\end{notation}

\begin{lemma}
  \label{thm:section}
  Let $E\to\bbS^1$ be a analytic line bundle such
  that (1) $E^\ast = E$, and
  (2) for each $\lambda\in\bbS^1$ except possibly at finitely many points,
  there exists $Y\in E(\lambda)$ which is positive definite.  Then there
  exists a analytic section $X$ of $E$ such that $X=X^\ast$,
  $X$ is positive semidefinite on $\bbS^1$, and $\det X\notequiv 0$.
\end{lemma}

\begin{proof}
  Let $X_1$ be a nowhere vanishing section of $E$.
  Then there exists $\alpha\in\bbC^\ast$ such that
  $X_2 = \alpha X_1+(\alpha X_1)^\ast\notequiv 0$,
  and $X_2$ is a section of $E$ satisfying $X_2^\ast=X_2$.

  For any $\lambda\in\bbS^1$ at which there exists $Y\in E(\lambda)$
  which is positive definite,
  since $\dim E_\lambda=1$ and $Y\ne 0$,
  $X_2(\lambda)=c Y$ for some $c\in\bbC$.
  Since at $\lambda$, $X_2=X_2^\ast$ and $Y=Y^\ast$, $c\in\bbR$.
  Hence $X_2(\lambda)$ is either positive definite, negative definite
  or $0$ according as $c>0$, $c<0$ or $c=0$.

  Let $P=\{p_1,\dots, p_n\}\subset\bbS^1$ be the set of points
  at which $X_2$ switches between being positive and negative definite.
  Then $P$ is even.
  Let $f(\lambda)=\lambda^{-n}\prod_{j=1}^{2n}(\lambda-p_i)$.
  Let $p\in\bbS^1\setminus P$ be a point for which $X_2(p)$ is
  positive definite and let $g(\lambda)=f(\lambda)/f(p)$.
  Then $g$ is analytic, $g\notequiv 0$, $g^\ast=g$, and
  $X_2$ is positive or negative definite according as $g>0$ or $g<0$.
  Thus $X=gX_2$ satisfies $\det X\notequiv 0$ and $X=X^\ast$ and
  is positive definite except at $P$, and is hence is  positive semidefinite.
\end{proof}

\subsection{The gluing theorem}
\label{sec:glue1}

We prove the main unitarization result:
if a set of monodromies is unitarizable pointwise on $\bbS^1$, then it is
unitarizable by an $r$-dressing. In the context of DPW,
such a dressing closes the periods of the CMC immersion by
lemma~\ref{thm:close}.
The proof is based on lemmas~\ref{thm:birkhoff1}--\ref{thm:section}.

\begin{theorem}\label{thm:glue}
  Let $M_k:\bbS^1\to\matGL_2(\bbC)$ $(k\in\{1,\dots,n\})$
  be analytic maps such that the set
  $\{M_1,\dots,M_n\}$
  is nondegenerate and
  simultaneously unitarizable pointwise on $\bbS^1$
  except possibly at a finite subset of $\bbS^1$.
  Then there exists an analytic map $C\in\LoopprX$
  for which $CM_kC^{-1}$ extends analytically across
  $\{\det C= 0\}$ and is in $\LoopuGL{1}$.

  Moreover, $C$ is unique up to multiplication by a scalar function
  $\bbS^1\to\bbC$ which is the boundary of an analytic function
  $\calD_1\to\bbC^\ast$.
\end{theorem}

\begin{proof}
  Let $L_\lambda:\mattwo(\bbC)\to (\mattwo(\bbC))^n$ be
  the linear map defined by
  \begin{equation*}
    L_\lambda(X)=(XM_1-{M_1^\ast}^{-1}X,\dots,XM_n-{M_n^\ast}^{-1}X).
  \end{equation*}
  $L_\lambda$ depends analytically on $\lambda\in\bbC^\ast$ because
  $M_j$ do.  $L_\lambda$ is constructed so its kernel is the
  ``square'' of a unitarizer in the following sense: an analytic map
  $C:\bbS^\ast\to\matGL_2(\bbC)$ satisfies $C^\ast C\in\ker
  L_\lambda$ if and only if $CM_jC^{-1}$, $j\in\{1,\dots,n\}$ satisfy
  the reality condition~\eqref{eq:reality}.
  
  By lemma~\ref{thm:dim1}, for $\lambda\in\bbS^1$ for which
  $\{M_1,\dots,M_n\}$
  is nondegenerate, $\dim\ker L_\lambda=1$.  By lemma~\ref{thm:mindim}(i),
  there exists a trivial analytic line bundle $E\to\bbS^1$
  such that $E_\lambda=\ker L_\lambda$ except possibly at a finite subset
  of $\bbS^1$, where $E_\lambda\subset\ker L_\lambda$. $E$
  satisfies conditions (1) and (2) in the hypothesis of
  lemma~\ref{thm:section}, so by that theorem, there exists a
  analytic section $X$ of $E$ with the properties $X=X^\ast$,
  $X$ is positive semidefinite on $\bbS^1$, and
  $\det X\notequiv 0$.
  
  By lemma~\ref{thm:birkhoff2}, there exist a ``square root'' of $X$
  in the sense that there exist analytic maps
  $C\in\LoopprX$ and $f:\bbS^1\to\bbC$ such that
  $f X = C^\ast C$.
  Then $CM_j C^{-1}$
  satisfies the conditions of
  lemma~\ref{thm:extend-unitary}, so by that lemma
  it extends analytically across $\{\det C=0\}$ and
  is in $\LoopuGL{1}$.

  To show uniqueness,
  let $C_1,\,C_2$ be two such maps, and
  let $A=C_2C_1^{-1}$ Then $A\in\LooppGL{r}$ for every $r\in(0,\,1)$.
  For each $\lambda\in\bbS^1$ except possibly at a finite set
  $S\subset\bbS^1$,
  $A(\lambda)$ unitarizes the unitary matrices $C_1M_kC_1^{-1}|_\lambda$.
  By lemma~\ref{thm:isu}, $A(\lambda)\in\bbR^+\times\matU_2$.
  By lemma~\ref{thm:extend-unitary},
  $A|_{\bbS^1}=fU$ for some meromorphic
  function $f:\bbS^1\to\bbR$ and analytic $U\in\LoopuGL{1}$.
  \Note{Check this. Why does $f$ extend to $\calD_1$?}
  For some $r$ close to $1$, the $r$-Iwasawa factorization
  of $A|_{\calC_r}$ is then $A|_{\calC_r}=U\cdot(f\id)$.
  But $A|_{\calC_r}\in\LoopprGL{r}$,
  so $U=\id$ and $C_2=fC_1$.
\end{proof}

\begin{figure}[ht]
  \centering
  \includegraphics[width=175pt]{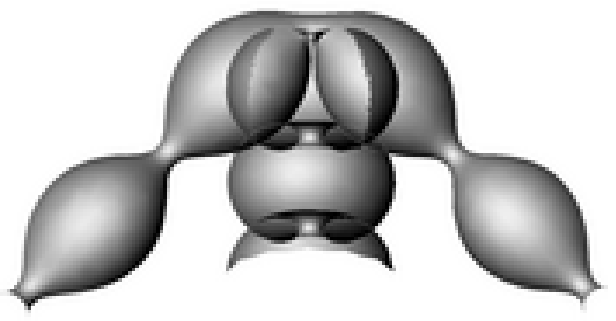}
  \includegraphics[width=175pt]{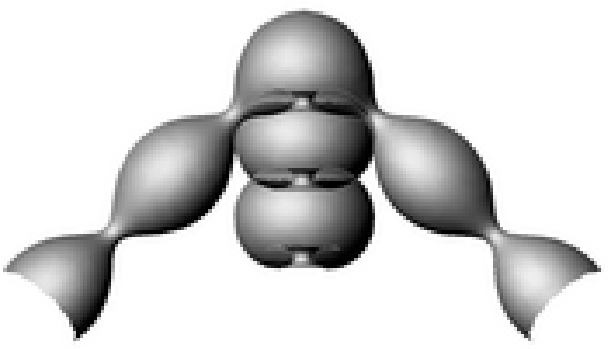}
\caption{A pair of CMC trinoids with one nodoid end (necksizes
$\bigl(\tfrac{1}{6},\,\tfrac{1}{6},\,-\tfrac{1}{6}\bigr)$).
The unduloid ends can be thought of as pulling outward along their axes, while
the nodoid end pushes upward, in static equilibrium.}
\end{figure}

\typeout{== trinoids.tex =============================================}
\section{Constructing trinoids}
\label{sec:trinoid}

Constructing trinoids is in the following steps:

1. Write down a family of DPW potentials on the thrice-punctured
sphere which are locally gauge-equivalent to perturbations of the
Delaunay DPW potential at each puncture
(definition~\ref{def:trinoid}).

2. Show that the monodromy representation is unitarizable pointwise
for $\lambda\in\bbS^1$
(theorems~\ref{thm:unitary3} and~\ref{thm:triunitary}).

3. Construct by the gluing theorem~\ref{thm:glue}
a dressing for which the monodromy representation is unitary on $\bbS^1$.
This dressing will close the three ends of the surface.

4. Show by the asymptotics theorem~\ref{thm:asymptotic} that
the three ends are asymptotically Delaunay.

\subsection{Trinoid potentials}
\label{sec:potential}

In this section a family of potentials is defined which will
be used produce trinoids via the DPW construction.
Near the punctures
the potentials are local perturbations of Delaunay potentials
via gauge equivalence.
The family is parametrized by the three asymptotic Delaunay weights
and has four connected components,
divided according as the necksizes are positive or negative:
\Wppp, \Wppm, \Wpmm, \Wmmm.

\begin{definition}
  \label{def:trinoid}
  Let $\Sigma=\bbP^1\setminus\{0,\,1,\,\infty\}$.
  Let $w_1,\,w_2,\,w_3\in(-\infty,\,1)\setminus\{0\}$
  and $W=(w_0,\,w_2,\,w_3)$.
  Let $n_j=\half(1-\sqrt{1-w_j})$, $j=1,\,2,\,3$ and suppose that
  $n_j$ and $w_j$ satisfy the inequalities
  \begin{align}
    \begin{split}
      \label{eq:neck}
      & \abs{n_1}+\abs{n_2}+\abs{n_3} \le 1\\
      & \abs{n_i} \le \abs{n_j}+\abs{n_k},\quad\{i,\,j,\,k\}=\{1,\,2,\,3\}\\
    \end{split}\\
    \label{eq:weight}
    & \abs{w_i} \le \abs{w_j} + \abs{w_k},
    \quad\{i,\,j,\,k\}=\{1,\,2,\,3\}.
  \end{align}
  Define $\xi_W\in\Potentialsl{\Sigma}{1}$ by
  \begin{equation}
    \xi_W=\begin{pmatrix}
      0 & \lambda^{-1}dz \\ (\lambda-1)^2Q_W/dz & 0
    \end{pmatrix}
  \end{equation}
  where
  \begin{equation}
    Q_W =\frac{w_3z^2-(w_1-w_2+w_3)z+w_1}{16z^2(z-1)^2}dz^2
  \end{equation}
is the unique meromorphic quadratic differential on
$\bbP^1$ whose only poles are double poles at $0,\,1,\,\infty$
with respective quadratic residues $w_k/16$.
By remark~\ref{rem:dpw}, the Hopf differential of the resulting CMC immersion
will be $-2H^{-1}Q_W\lambda^{-1}$.
\end{definition}

\subsection{Local gauge}

We show that the double pole of a trinoid potential can be gauged
to a simple pole with Delaunay residue and, after a coordinate
change, no constant term.

\begin{lemma}
  \label{thm:gauge-trinoid}
  Let $\xi_W\in\mathcal{T}$ be a trinoid potential.
  Then for each end $p\in\{0,\,1,\,\infty\}$
  there exists
  a neighborhood $U$ of $p$,
  an analytic map $g:U^\ast\to\LooppGL{1}$ and a
  conformal coordinate $\tilde{z}:U\to\bbC$ with $\tilde{z}(p)=0$,
  such that the expansion of $\gauge{\xi_W}{g}$ is
  \begin{equation}
    \label{eq:gauge-trinoid}
    \begin{pmatrix}
      0 & a\lambda^{-1}+b \\ b+a\lambda & 0
    \end{pmatrix}\frac{d\tilde{z}}{\tilde{z}}+O(\tilde{z})d\tilde{z}.
  \end{equation}
\end{lemma}

\begin{proof}
  Let $\mu_{w_1}$ as in equation~\eqref{eq:mu}.
  There exists $a,\,b\in\mathbb{R}$ with $|a|\ge|b|$ satisfying
  $pp^\ast=\mu^2$, where $p=a\lambda^{-1}+b$.
  Let
  \begin{equation}
    \label{eq:g}
    g_1=
    \begin{pmatrix}
      z^{1/2} & 0 \\ 0 & z^{-1/2}
    \end{pmatrix},\quad
    g_2=
    \begin{pmatrix}
      1 & 0 \\ -\tfrac{1}{2}\lambda & \lambda p
    \end{pmatrix}.
  \end{equation}
  Then $\gauge{\xi_W}{g_1g_2}$ has a simple pole at $z=0$ and residue
  as in equation~\eqref{eq:gauge-trinoid}.
  Let
  \begin{equation*}
    g_3=\id+\frac{k}{2}
    \begin{pmatrix}
      -1 & 0 \\ p^{-1} & 1
    \end{pmatrix}z,\quad
    k=\frac{w_1+w_2-w_3}{2w_1},\quad
    z = \tilde{z} - k \tilde{z}^2,
  \end{equation*}
  \Note{check this for general choice of ends.}
  be the gauge and coordinate change constructed by
  lemma~\ref{thm:gaugeconst}.
  Then $g=g_1g_2g_3$ and $\widetilde{z}$
  are the required gauge and coordinate change.
\end{proof}

\subsection{Gauge-equivalent trinoid potentials}

We present two gauge-equivalent forms of the trinoid potentials of
definition~\ref{def:trinoid}.
Lemma~\ref{thm:lowertriangulargauge} shows that any potential may be gauged
to an off-diagonal form with a prescribed upper-right entry.

\begin{lemma}
  \label{thm:lowertriangulargauge}
  Let $\Sigma$ be Riemann surface and $\widetilde{\Sigma}$ its
  universal cover. Let $r\in(0,\,1]$ and let
  $\xi\in\Potentialsl{\Sigma}{r}$ be given by
  \[
  \xi =
  \begin{pmatrix}
    c & \lambda^{-1}a \\ b & -c
  \end{pmatrix}\omega.
  \]
  where $a,\,b,\,c$ are mermorophic functions on $\Sigma$ depending on
  $\lambda$ and $\omega$ is a $\lambda$-independent
  meromorphic 1-form on $\Sigma$.
  Let $s\in(0,\,r]$ such that $a$ has no zeros in
  $\{0\le\lambda\le s\}$. Then the map
  $g:\widetilde{\Sigma}\to\LooppSL{s}$ defined by
  \[
  g = 
  \begin{pmatrix}
    a^{1/2} & 0 \\
    \lambda\left( \frac{d(a^{-1/2})}{\omega} - c a^{-1/2} \right) & a^{-1/2}
  \end{pmatrix}
  \]
  gauges $\xi$ to
  \[
  \gauge{\xi}{g}=
  \begin{pmatrix}
    0 & \lambda^{-1}\omega \\ Q/\omega & 0
  \end{pmatrix}
  \]
  for some meromorphic quadratic differential $Q$ on $\Sigma$.
\end{lemma}

\begin{lemma}
  \label{thm:gaugetodel}
  $\xi_W\in\calT$ can be gauged globally to Fuchsian system
  with hermitian residues as in~\cite{schmitt1}.
  This gauge introduces extra poles with weight $0$ and monodromy $-\id$.
\end{lemma}

\begin{proof}
  We provide the gauge in the case of three positive weights.
  The proof in the other cases is similar.

  Potentials in the family in~\cite{schmitt1} are of the form
  \begin{equation*}
    \xi = \begin{pmatrix}
      \gamma & \alpha\lambda^{-1}+\beta\\
      \beta + \alpha\lambda & -\gamma
    \end{pmatrix}
  \end{equation*}
  where
  $W=(w_1,\,w_2,\,w_3)\in\bbR^3$,
  \begin{equation*}
    \begin{split}
      w &= \half(w_1+w_2+w_3),\\
      r_k &= \frac{\sqrt{w-w_i}\,\sqrt{w-w_j}}{4\sqrt{w-w_k}},\quad
      \{i,\,j,\,k\}=\{1,\,2,\,3\},\\
      r &= r_1+r_2+r_3,\\
      p &= -\frac{1}{2r}+\sqrt{\frac{1}{4r^2}-1},
    \end{split}
  \end{equation*}
  taking positive square roots, and
  \begin{equation*}
    \begin{split}
      \alpha &= a\, dz = \left(\frac{r_1}{z}+\frac{r_2}{z-1}\right)dz\\
      \beta &= b\, dz = \left(\frac{r-r_1}{z}+
        \frac{r-r_2}{z-1}-\frac{r}{z-\frac{r_1}{r_1+r_2}}\right)dz\\
      \gamma &= \half(p-p^{-1})(\alpha+\beta)
    \end{split}.
  \end{equation*}
  The potential $\xi$ has simple poles at
  $(0,\,1,\,\infty,\,\frac{r_1}{r_1+r_2})$
  with residues of the form~\eqref{eq:del} with respective
  weights $(w_1,\,w_2,\,w_3,\,0)$.
  Let
  \[
  h = \frac{1}{\sqrt{1-\lambda}}
  \begin{pmatrix}
    1 & p \\ p^{-1}\lambda & 1
  \end{pmatrix}.
  \]
  and $g$ be the gauge of lemma~\ref{thm:lowertriangulargauge}
  obtained from $\gauge{\xi}{h}$, taking $\omega=dz$ in that lemma.
  Then $\gauge{\xi}{hg}\in\calT_W$.
\end{proof}

\begin{lemma}
  \label{thm:gaugetodorf}
  The family of trinoid potentials in~\cite{dw-trinoid} is gauge equivalent
  to the family $\calT$.
\end{lemma}

SHOW HOW TO GAUGE TO HYPERGEOMETRIC EQUATION INSTEAD.

\begin{proof}
  Potentials in the family in~\cite{dw-trinoid} are of the form
  \begin{equation*}
    \xi =
    \begin{pmatrix}
      0 & \sigma \\ \tau & 0
    \end{pmatrix},
  \end{equation*}
  where
  $W=(w_0,\,w_1,\,w_\infty)\in\bbR^3$,
  $a_0,\,a_1\in\bbZ$, and $\omega$ is an analytic loop on $\bbS^1$ which
  extends to a holomorphic function on $\calD_1^\ast$ with no zeros,
  and extends meromorphically to $0$ with
  $\ord_0\omega=-1$,
  \begin{equation*}
    \begin{split}
      \sigma &= \omega z^{-a_0}(z-1)^{-a_1}\\
      -\sigma\tau &=
      \frac{b_0}{z^2}+\frac{b_1}{(z-1)^2}+\frac{c}{z}-\frac{c}{z-1}
    \end{split}
  \end{equation*}
  and
  \begin{equation*}
    \begin{split}
      b_k &= ((a_k-1)/2)^2-\mu_k^2,\quad k=0,\,1\\
      c &= 1/4 - a_0 a_1/2 - \mu_0^2 - \mu_1^2 + \mu_\infty^2\\
      \mu_k &= \frac{1}{2}\sqrt{1+\frac{w_k(\lambda-1)^2}{4\lambda}},\quad
      k\in\{0,\,1,\,\infty\}.
    \end{split}
  \end{equation*}
  Let $g$ as in lemma~\ref{thm:lowertriangulargauge}
  taking $\omega=dz$ in that lemma.
  Then $\gauge{\xi}{g}\in\calT_w$.
\end{proof}

\subsection{Unitary monodromy on the thrice-punctured sphere}
\label{sec:unitarize}

In this section it is shown that given $M_1,\,M_2,\,M_3\in\matSL_2(\bbC)$
whose product is $\id$, the spherical triangle inequalities on the logs of
their eigenvalues are necessary sufficient for the simultaneous
unitarizability of $M_1,\,M_2,\,M_3$.
An equivalent condition in terms of the traces of the matrices
is given in~\cite{Goldman1988}.
Such inequalities are discussed in the context of holomorphic vector bundles
in~\cite{biswas}.

For a set of more than three matrices whose product is $\id$,
the spherical $n$-gon inequalities are necessary but not sufficient
conditions for simultaneous unitarizability. The case $n=3$ is
special in that the dimension of the set of conjugacy classes for
$M_1,\,M_2,\,M_3$ is the same as that of the eigenvalues.

\begin{lemma}[Spherical triangle inequalities]
  \label{thm:spheretri}
  Given $(\nu_1,\,\nu_2,\,\nu_3)\in(0,\,\half)^3$, there exists
  a nondegenerate spherical triangle on $\bbS^1$ with sides $2\pi\nu_k$ iff
  $(\nu_1,\,\nu_2,\,\nu_3)$ satisfy the spherical triangle inequalities
  \begin{equation}
    \label{eq:spheretri}
    \begin{split}
      & \nu_1+\nu_2+\nu_3 < 1,\\
      & \nu_i < \nu_j+\nu_k,\quad\{i,\,j,\,k\}=\{1,\,2,\,3\}.
    \end{split}
  \end{equation}
\end{lemma}

\begin{lemma}\label{thm:ulemma}\mbox{}
  (i) $M\in\matSL_2(\bbC)$ is unitarizable (notation~\ref{not:unit}) iff
  $\half\tr M\in(-1,\,1)$ or $M\in\{\pm\id\}$.

  (ii) Any $M\in\matSU_2$ can be written $M=\cos(2\pi\nu)+\sin(2\pi\nu)A$ with
  $\nu\in[0,\,\half]$ and $A\in\matsu_2$ with $\det A=1$.
\end{lemma}

\begin{theorem}
  \label{thm:unitary3}
  Let $M_1,\,M_2,\,M_3\in\matSL_2(\bbC)$
  with $M_1M_2M_3=\id$ and
  with eigenvalues $\exp(\pm 2\pi i\nu_k)$, $\nu_k\in(0,\,\half)$.
  Then $M_1,\,M_2,\,M_3$ are nondegenerate and
  simultaneously unitarizable iff
  the spherical triangle inequalities~\eqref{eq:spheretri} hold.
\end{theorem}

\begin{proof}
  Suppose $M_k$ are nondegenerate and simultaneously unitarizable,
  and let $C$ be a unitarizer,
  so that $CM_kC^{-1}\in\matSU_2$.
  Write $CM_kC^{-1} = x_k\id + y_k A_k$ as in lemma~\ref{thm:ulemma}(ii).
  The nondegeneracy assumption means the $A_k$ span $\matsu_2$.
  Identifying $\matsu_2\equiv\bbR^3$,
  let $P_k$ be the planes perpendicular to
  $A_k$ through $0$.
  The planes intersect $\bbS^2$
  forming eight spherical triangles;
  consider one of the spherical triangle $\Delta$ with
  side lengths less than $\pi$.
  An spherical trigonometry argument shows that the side lengths of
  $\Delta$ are $\nu_1,\,\nu_2,\,\nu_3$, so by lemma~\ref{thm:spheretri},
  the spherical triangle inequalities~\eqref{eq:spheretri} hold.

  Conversely, given
  $(\nu_1,\,\nu_2,\,\nu_3)\in(0,\,\half)^3$ satisfying
  the spherical triangle inequalities, by lemma~\ref{thm:spheretri}
  there exists a nondegenerate spherical triangle on $\bbS^2$
  with side lengths $\nu_1,\,\nu_2,\,\nu_3$. Let
  $A_k$ be the normals to the planes through the sides.
  Then $M_k=\cos(2\pi i\nu_k)\id + \sin(2\pi i\nu_k) A_k$
  are nondegenerate and unitary, and
  a spherical trigonometry argument shows that $M_1M_2M_3=\id$.

  It remains to show that a choice $(\nu_1,\,\nu_2,\nu_3)$ determines
  a unique conjugacy class of $(M_1,\,M_2,\,M_3)$.
  If $(N_1,\,N_2,\,N_3)$ is another triple with the same traces,
  we can assume by conjugation that $N_1=M_1$, and need to show that
  $N_2$ is conjugate to $M_2$ by a commutator of $M_1$.
  A computation shows
  \begin{equation*}
    M_1^\circ M_2^\circ + M_2^\circ M_1^\circ = 2(t_3-t_1t_2)\id,
  \end{equation*}
  where $X^\circ$ denotes $\tracefree(X)$.
  Since $\nu_1\notin\{0,\,\half\}$, then $M_1\notin\{\pm\id\}$
  and $M_1^\circ\ne 0$. Fixing $M_1^\circ$, the equation is linear
  in $M_2^\circ$ and has a 2 complex dimensional solution space.
  \Note{Assuming tracefree but not det=1}
  Since if $M_2^\circ$ is a solution, so is $CM_2^\circ C^{-1}$ for any
  commutator $C$ of $M_1$, and the set of such commutators is
  also a 2 complex dimensional linear space, the set of solutions is a
  single orbit under conjugation by commutators of $M_1$. The result follows.
\end{proof}

\subsection{Unitarization of trinoid monodromy pointwise on $\bbS^1$}
\label{sec:moduli}

We compute the eigenvalues of the monodromy for a potential $\xi_W\in\calT$.

\begin{lemma}
  \label{thm:trinoid-mono}
  Let $\xi_W\in\calT$ be a trinoid potential,
  $\Phi$ a solution to the ODE $d\Phi=\Phi\xi_W$,
  and $M_1,\,M_2,\,M_3$ the monodromy of $\Phi$ at $0,\,1,\,\infty$
  respectively.
  Then the eigenvalues of $M_k$ are $\exp(\pm2\pi\nu_{w_k})$,
  where $\nu_w$ is defined by
  \begin{equation}
    \label{eq:nu}
    \nu_w = \frac{1}{2}-\frac{1}{2}\sqrt{1+\frac{w(\lambda-1)^2}{4\lambda}}.
  \end{equation}
\end{lemma}

\begin{proof}
  By lemma~\ref{thm:gauge-trinoid}, $\xi_W$ is locally gauge-equivalent
  to a potential $\eta$ of the form of equation~\eqref{eq:gauge-trinoid}.
  Let $M_{\xi_W}$ and $M_\eta$ be the respective monodromy representations
  of $\xi_W$ and $\eta$.
  By lemma~\ref{thm:eigen}(i), the eigenvalues of the monodromy of
  $\eta$ are $\exp(\pm2\pi i (\half-\nu_w))$, where
  $\nu_w$ is given by equation~\eqref{eq:nu}.
  By lemma~\ref{thm:gauge}(ii), $M_{\xi_W}=-M_\eta$,
  hence the eigenvalues of the monodromy representation of $\xi_W$
  are $\exp(\pm 2\pi i\nu_w)$.
\end{proof}

Necessary
and sufficient conditions are found that a monodromy representation with these
eigenvalues be unitarizable for every $\lambda\in\bbS^1$
(conditions~\eqref{eq:neck}--\eqref{eq:weight}).
The inequalities on the necks $n_i$ are
the spherical triangle inequalities on the eigenvalues evaluated
at $\lambda=-1$. The inequalities on the weights $w_i$ are
implied by the balancing formula, according to which the sum of the
end forces (the end axes in $\matsu_2$ with length $w_k$) is $0$.

\begin{notation}\label{not:T}
  Let $T_0\subset\bbR^3$ be the bounded set with tetrahedral boundary defined
  by
  \begin{equation*}
    \begin{split}
      & \nu_1+\nu_2+\nu_3 \le 1,\\
      & \nu_i \le \nu_j+\nu_k,\quad\{i,\,j,\,k\}=\{1,\,2,\,3\}
    \end{split}
  \end{equation*}
  and let $T$ be the orbit of $T_0$ by the action of the group
  generated by the transformations $\nu_k\mapsto \nu_k + 1$ and
  $\nu_k\mapsto -\nu_k$.  
\end{notation}

\begin{lemma} 
  \label{thm:intet}
  Let $\nu_{w_k}$ be defined by equation~\eqref{eq:nu}
  and $\nu=(\nu_1,\,\nu_2,\,\nu_3)$
  Then $\nu\in T$ for all $\lambda\in\bbS^1$ iff
  the inequalities~\eqref{eq:neck} and~\eqref{eq:weight} are satisfied.
\end{lemma}

\begin{proof}
  Assume equations~\eqref{eq:neck} and~\eqref{eq:weight} are satisfied.
  Define
  \begin{equation*}
    \label{eq:f}
    \begin{split}
      \rho_k &= \half-\half\sqrt{1-w_k x},\quad \{i,\,j,\,k\}=\{1,\,2,\,3\}\\
      f &= \abs{\rho_1}+\abs{\rho_2}+\abs{\rho_3}\\
      f_i &= -\abs{\rho_i}+\abs{\rho_j}+\abs{\rho_k},\quad
      \{i,\,j,\,k\}=\{1,\,2,\,3\}.
    \end{split}
  \end{equation*}
  The terms in $f$ are increasing, so $f$ is increasing, so
  $n_1+n_2+n_3\le 1$ implies that $f\le 1$ on $[0,\,1]$.
  Hence $\nu_1+\nu_2+\nu_3\le 1$ on $\bbS^1$.

  In the case $0<w_1\le w_2$ or $w_2\le w_1<0$, $f_1$ is increasing,
  so $n_1\le n_2+n_3$ implies $f_1$ is non-negative on $[0,\,1]$.
  Hence $\nu_1\le\nu_2+\nu_3$ on $\bbS^1$.

  We require the following fact: the function
  $\rho_2/\rho_1$ extends to a $C^\infty$ function
  at $0$, and, if $w_2>w_1$, then $\abs{\rho_2/\rho_1}$
  is strictly increasing.

  In the case $w_1\ge w_2$, $w_1\ge w_2$, the above fact implies that
  that $f_1/\abs{\rho_1}$ is non-increasing.
  $n_1\le n_2+n_3$ implies that $f_1/\abs{\rho_1}$ is non-negative at
  $1$, so $f_1/\abs{\rho_1}$, and hence $f_1$, is non-negative on $[0,\,1]$.
  Hence $\nu_1\le\nu_2+\nu_3$ on $\bbS^1$.

  In the case $w_1\le w_2$, $w_1\le w_2$, the above fact implies that
  that $f_1/\abs{\rho_1}$ is non-decreasing.
  But $(f_1/\abs{\rho_1})(0)=-1+\abs{w_2/w_1}+\abs{w_3/w_1} \ge 0$, so
  $f_1/\abs{\rho_1}$ is non-negative on $[0,\,1]$.
  Hence $\nu_1\le\nu_2+\nu_3$ on $\bbS^1$.

  Symmetric arguments for the other cases imply that
  $\nu\in T$.

  The proof of the converse is omitted.
\end{proof}

\begin{lemma}
  \label{thm:w3}
  If the conditions~\eqref{eq:neck}--\eqref{eq:weight}
  are satisfied, then $w_k > -3$, $k=1,\,2,\,3$.
\end{lemma}

\begin{proof}
  The inequalities $\abs{n_i}\le\abs{n_j}+\abs{n_k}\le 1-\abs{n_i}$
  imply $\abs{n_i}\le\half$.
  Hence $w_i\ge -3$.
  Suppose $w_3= -3$, so $n_3=-\half$.
  By the above inequalities, $\abs{n_1}+\abs{n_2}=\half$.
  Using $w_r=4n_r(1-n_r)$, the inequality
  $\abs{w_3}\le\abs{w_1}+\abs{w_2}$ implies
  $\tfrac{1}{4}\le -n_1\abs{n_1}-n_2\abs{n_2}$.
  An examination of cases according to the signs of $n_1$, $n_2$
  shows that this is satisfied only if $n_1=0$ or $n_2=0$.
\end{proof}

The following theorem, the main theorem of the section, shows that
the monodromy representation of a trinoid potential is pointwise
unitarizable on $\bbS^1$.

\begin{theorem}
  \label{thm:triunitary}
  Let $\xi_W\in\calT$, and let $\Phi$ a solution to the ODE $d\Phi=\Phi\xi_W$
  such that $\Phi(p)$ is holomorphic on $\calA_0$ for some $p$ in the
  universal cover of $\Sigma$.
  Then the monodromy representation of $\Phi$ is nondegenerate and pointwise
  unitarizable on $\bbS^1$ except possibly at finitely many points.

  Conversely, the conditions~\eqref{eq:neck}--\eqref{eq:weight}
  are necessary in order for the
  monodromy representation of $\Phi$ to be nondegenerate and pointwise
  unitarizable on $\bbS^1$ except possibly at finitely many points.
\end{theorem}

\begin{proof}
  By lemma~\ref{thm:trinoid-mono}, 
  the eigenvalues of $M_k$ are $\exp(\pm2\pi\nu_{w_k})$,
  where $\nu_{w_k}$ are defined by equation~\eqref{eq:nu}.

  A necessary condition for the degeneracy of $\{M_k\}$ on $\bbS^1$ is
  $\nu\in\del T$, but this occurs only at finitely many points on
  $\bbS^1$.

  By the definition of $\calT$, $w_1,\,w_2,\,w_3$ satisfy the
  neck and weight inequalities~\eqref{eq:neck}--\eqref{eq:weight}.

  Let $S=\{\lambda\in\bbS^1\suchthat \nu\in\del T\}$. Then $S$ is finite.
  Then the following are equivalent:
  (1) $M_1,\,M_2,\,M_3$ are irreducible and simultaneously unitarizable on
  $\bbS^1\setminus S$.
  (2) $(\nu_1,\,\nu_2,\,\nu_3)\in\interior{T}$ for all
  $\lambda\in\bbS^1\setminus S$
  (theorem~\ref{thm:unitary3}).
  (3) $(\nu_1,\,\nu_2,\,\nu_3)\in T$ for all $\lambda\in\bbS^1$.
  (4) The inequalities~\eqref{eq:neck} and~\eqref{eq:weight} hold
  (lemma~\ref{thm:intet}).
\end{proof}

\subsection{Main theorem}
\label{sec:main}

\begin{theorem}\label{thm:main}
  Let $\Sigma=\bbP^1\setminus\{0,\,1,\,\infty\}$,
  let $w_1,\,w_2,\,w_3\in(-\infty,\,1]\setminus\{0\}$
  and $n_j=\half(1-\sqrt{1-w_j})$, $j=1,\,2,\,3$ and assume
  \begin{align*}
    \begin{split}
      & \abs{n_1}+\abs{n_2}+\abs{n_3} \le 1\\
      & \abs{n_i} \le \abs{n_j}+\abs{n_k},\quad\{i,\,j,\,k\}=\{1,\,2,\,3\}\\
    \end{split}\\
    & \abs{w_i} \le \abs{w_j} + \abs{w_k},
    \quad\{i,\,j,\,k\}=\{1,\,2,\,3\}.
  \end{align*}
  Then there exists a conformal CMC immersion $f:\Sigma\to\bbR^3$
  with three ends which are asymptotic to half Delaunay surfaces
  with weights $w_1,\,w_2,\,w_3$.
\end{theorem}

\begin{proof}
  Let $W=(w_1,\,w_2,\,w_3)$ and let
  $\xi_W\in\calT$ be a trinoid potential (definition~\ref{def:trinoid}).
  Let $\widetilde{\Sigma}\to\Sigma$ be the universal cover of $\Sigma$
  and $\Gamma$ the group of deck transformations for this cover.
  let $\Phi\in\LoopGL{1}$ a nonsingular solution to the ODE $d\Phi=\Phi\xi_W$
  which extends analytically to $\calA_0$.
  Let $M_1,\,M_2,\,M_3$ the monodromies of $\Phi$ at $0,\,1,\,\infty$
  respectively.

  \emph{Step 1: Closing the ends.}
  By theorem~\ref{thm:triunitary},
  the set $M_1,\,M_2,\,M_3$ is nondegenerate and pointwise
  simultaneously unitarizable on $\bbS^1$ except possibly at a finite
  subset of $\bbS^1$.
  Thus by the gluing theorem~\ref{thm:glue} there exists an analytic map
  $C\in\LoopprX$ for which $CM_kC^{-1}$ extends analytically across
  $\{\det C= 0\}$ and is in $\LoopuGL{1}$.
  Let $r=(0,\,1)$ and
  $f_\lambda = \Sym{\lambda}{\Uni{r}{C\Phi}}$.
  By lemma~\ref{thm:rindependent}, $f$ is independent of the choice of $r$.

  By lemma~\ref{thm:eigen}(ii), the eigenvalues of $CM_kC^{-1}$
  on $\bbS^1$ are $\exp(\pm2\pi\nu_{w_k})$. These by construction
  satisfy equation~\eqref{eq:eigenclose}, so
  by lemma~\ref{thm:monoeigen}, $CM_kC^{-1}$
  satisfy the closing conditions~\eqref{eq:close2}.
  Hence by theorem~\ref{thm:close}, $f_1$ is closed in the sense that
  $\tau^{*}f_1=f_1$ for all $\tau\in\Gamma$.

  \emph{Step 2: Delaunay asymptotics.}
  Choose an end $p\in\{0,\,1,\,\infty\}$.
  By lemma~\ref{thm:w3}, the corresponding weight $w_k$ satisfies $w_k>-3$.
  By lemma~\ref{thm:gauge-trinoid}, there exists a gauge $g$ in a
  punctured neighborhood of $p$
  such that after a coordinate change, $\gauge{\xi_W}{g}$ has no constant
  term in its series expansion.
  Since $M_{C\Phi}\in\LoopuGL{r}$ by the construction of $C$,
  and $M_{C\Phi}=-M_{C\Phi g}$ by lemma~\ref{thm:gauge}(ii),
  then $M_{C\Phi g}\in\LoopuGL{r}$.
  Hence by the asymptotics theorem~\ref{thm:asymptotic}, 
  $\Sym{1}{\Uni{r}{C\Phi g}}$
  is asymptotic to half Delaunay
  surfaces at its ends $(0,\,1,\,\infty)$ with respective weights
  $w_1,\,w_2,\,w_3$.
  By lemma~\ref{thm:gauge}(iii) the same is true for $f_1$.
\end{proof}

The following theorem discusses the symmetry groups of the trinoids
constructed in theorem~\ref{thm:main}.

\begin{theorem}
  \label{thm:symmetry}
  (i) Each trinoid in the family
  constructed in theorem~\ref{thm:main} has a plane of reflective symmetry
  which fixes each end.
  (ii) Each isosceles trinoid in the family has a further plane of reflective
  symmetry perpendicular to this plane which exchanges the equal ends and
  fixes the third end.
  (iii) Each equilateral trinoid in the family has the order-12 symmetry
  group of an equilateral triangle slab.
\end{theorem}

The proof of this theorem will be found in~\cite{schmitt3}, which discusses
gauge symmetries in the general context of $n$-noids.

\begin{figure}[ht]
  \centering
  \includegraphics[width=175pt]{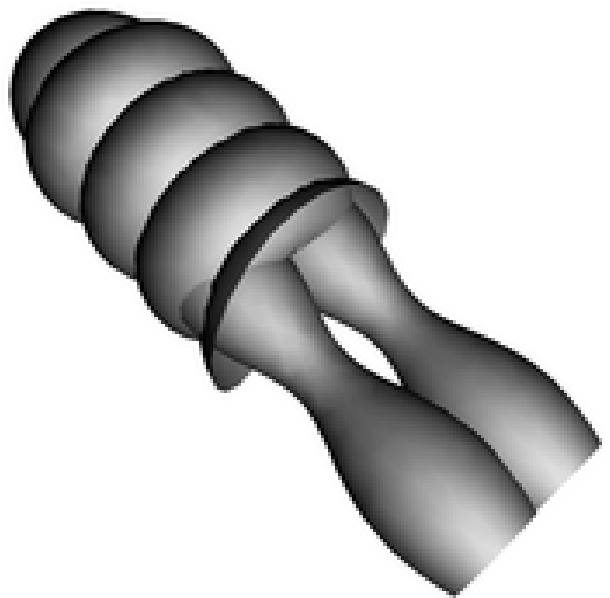}
  \includegraphics[width=175pt]{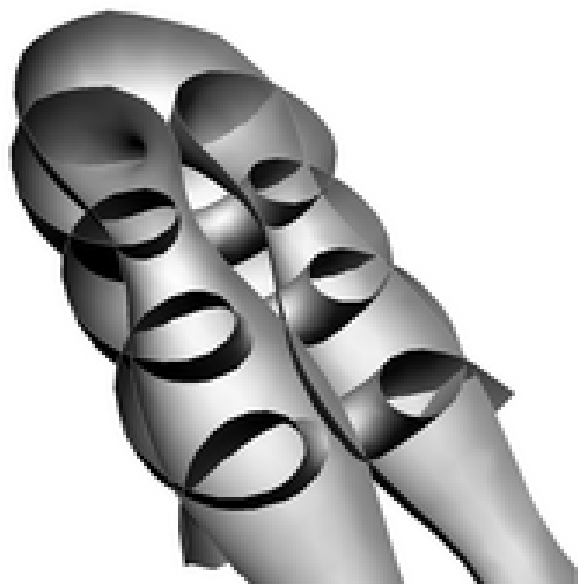}
\caption{Two views of a trinoid with one nodoid end (necksizes
$\bigl(\tfrac{1}{3},\,\tfrac{1}{3},\,-\tfrac{1}{3}\bigr)$).
The weight of the nodoid end is twice that of each unduloid end
but opposite in sign,
so the end axes are parallel as required by balancing.}
\end{figure}

\begin{figure}[ht]
  \centering
  \includegraphics[width=175pt]{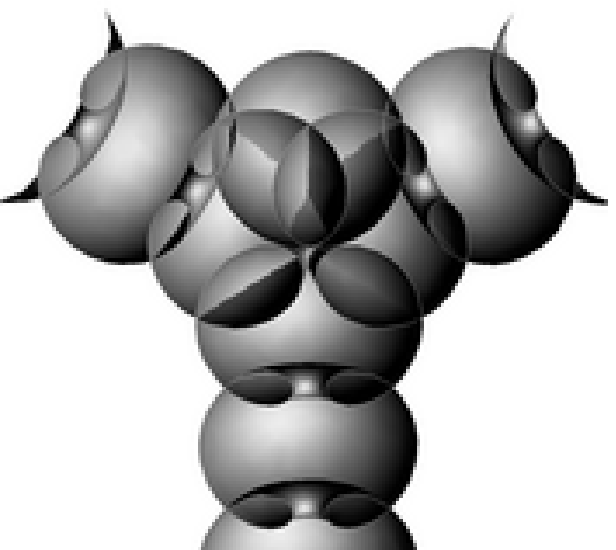}
  \includegraphics[width=175pt]{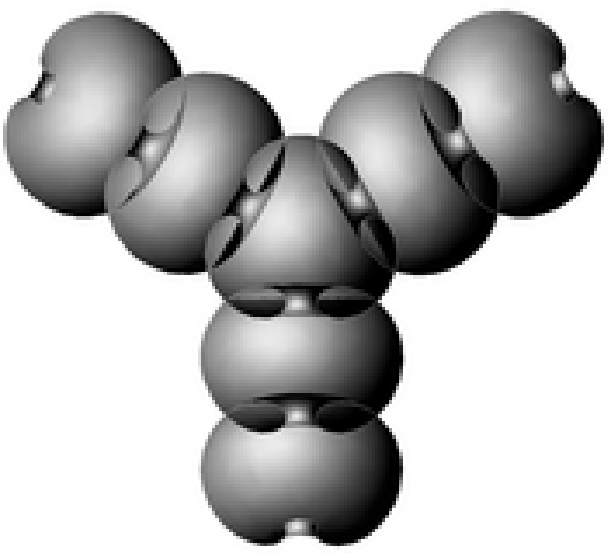}
\caption{A pair of CMC trinoids with three nodoid ends (necksizes
$\bigl(-\tfrac{1}{4},\,-\tfrac{1}{4},\,-\tfrac{1}{4}\bigr)$).}
\end{figure}

\typeout{== question.tex =============================================}
\section{Open questions}\label{sec:questions}

1. Computer experiments indicate that the trinoids in the subfamily
with embedded ends are Alexandrov embedded.

2. B\"ackland transformations can be applied to Delaunay surfaces
to obtain bubbletons~\cite{SW1}. Construct B\"ackland transformations of
CMC $n$-noids.

3. Classify the CMC trinoids.

4. Construct and classify the CMC $n$-noids. This will involve unitarizing
the monodromy representation on the $n$-punctured sphere.

5. Construct and classify $n$-noids with genus $>0$.

\bibliographystyle{amsplain}
\bibliography{ref}

\end{document}